\documentclass[a4paper,10pt]{article}
\usepackage{amssymb}
\newtheorem{Th}{Theorem}

\newtheorem{Co}{Corollary}
\newtheorem{Lm}{Lemma}

\newtheorem{Rm}{Remark}

\newcommand{\be}{\begin{equation}}
\newcommand{\ee}{\end{equation}}
\newcommand{\R}{\mathbb{R}}
\newcommand{\N}{\mathbb{N}}

\newcommand{\Z}{\mathbb{Z}}
\newcommand{\Q}{\mathbb{Q}}

\newcommand{\cqfd}
{%
\mbox{}%
\nolinebreak%
\hfill%
\rule{2mm}{2mm}%
\medbreak%
\par%
}
\newcommand\res{\mathop{\hbox{\vrule height 7pt width .5pt depth 0pt
\vrule height .5pt width 6pt depth 0pt}}\nolimits}

\def\lf{\left}
\def\rg{\right}

\def\ds{\displaystyle}

\def\p{\partial}

\def\res{\mathop{\hbox{\vrule height 7pt width .5pt 
depth 0pt\vrule height .5pt width 6pt depth 0pt}}\nolimits}
\newcommand{\HF}{\mathcal{H}}
\newcommand{\HH}{\mathbb{H}}
\newcommand{\B}{\mathbb{B}}
\newcommand{\M}{\mathbb{M}}
\newcommand{\sS}{\mathbb{S}}
\newcommand{\oO}{\mathbb{O}}
\newcommand{\pP}{\mathbb{P}}
\newcommand{\lseg}{\mathbf{[}\!\mathbf{[}}  
\newcommand{\rseg}{\mathbf{]}\!\mathbf{]}} 
\begin{document}
\title{ Sequential Weak Approximation for Maps of Finite Hessian Energy}

\author{Robert Hardt\footnote{Department of Mathematics, Rice University,
 Houston, TX 77251, USA. Research partially supported by the NSF.}\  and Tristan Rivi\`ere\footnote{Forschungsinstitut f\"ur Mathematik, ETH Zentrum,
CH-8093 Z\"urich, Switzerland.}}
\maketitle

{\bf Abstract :} {\it Consider the space $W^{2,2}(\Omega;N)$ of second order Sobolev mappings $\ v\ $ from a smooth domain $\Omega\subset\R^m$ to a compact Riemannian manifold $N$ whose Hessian energy   $\int_\Omega |\nabla^2 v|^2\, dx$  is finite.  Here we are interested in relations between the topology of $N$ and the $W^{2,2}$ strong or weak approximability of a $W^{2,2}$ map by a sequence of smooth maps from $\Omega$ to $N$. We treat in detail $W^{2,2}(\B^5,\sS^3)$ where we establish the \underline{sequential weak} $W^{2,2}$ density of $W^{2,2}(\B^5,\sS^3)\cap{\mathcal C}^\infty$.  The strong $W^{2,2}$ approximability of higher order Sobolev maps has been studied in the recent preprint \cite{BPV} of P. Bousquet, A. Ponce, and J. Van Schaftigen.  For an individual map $v\in W^{2,2}(\B^5,\sS^3)$, we define a number $L(v)$ which is approximately the total length required to connect the isolated singularities of a strong approximation $u$ of $v$ either to each other or to $\p\B^5$.  Then $L(v)=0$ if and only if  $v$ admits $W^{2,2}$ strongly approximable by smooth maps. Our critical result, obtained by constructing specific curves connecting the singularities of $u$, is the bound $\  L(u)\leq c\int_{\B^5}|\nabla^2 u|^2\, dx\ $. This allows us to construct, for the given Sobolev map $v\in W^{2,2}(\B^5,\sS^3)$, the desired $W^{2,2}$ weakly approximating sequence of smooth maps.
To find  suitable connecting curves for $u$, one uses the twisting of a $u$ pull-back normal framing of a suitable level surface of $u$}. 

\medskip

\noindent{\bf Math. Class.} 58D15, 46E35, 49Q99.

\smallskip

\tableofcontents


\section{Introduction}
To motivate our specific work on weak sequential approximability of  $W^{2,2}$ maps from $\B^5$ to $\sS^3$, we will first describe briefly the background and general problems.  Let $(N, g)$ be a compact Riemannian manifold. Via the Nash embedding theorem, one may assume that $N$ is a submanifold of some Euclidian space $\R^\ell$ and that the metric $g$ is induced by this inclusion.  One then has, for any open subset $\Omega$ of $\R^m$, $k\in\N$, and $p>1$, the nonlinear space of $k$th order, Sobolev maps
$$
W^{k,p}(\Omega,N)\ =\ \{u\in W^{k,p}(\Omega,\R^\ell )\ :\ u(x)\in N\ \rm{for\ almost\ every\ } x\in\Omega \},
$$
where $W^{k,p}(\Omega,\R^\ell )$ denotes the Banach space of $\R^\ell$-valued, order $k$, Sobolev functions on $\Omega$ with norm
$\|u\|_{W^{k,p}}\ =\ \left[\sum_{j=0}^k\left(\int_\Omega|\nabla^ju|^p\,dx\right)^{2/p}\right]^{1/2}$. 

\subsection{Strong Approximation} A basic question concerning the spaces $W^{k,p}(\Omega,N)$ is the {\it approximability of these maps by a sequence of smooth maps} of $\Omega$ into $N$. The issue involves the possible discontinuities in a Sobolev map because any continuous Sobolev map may be approximated strongly in the Sobolev norm. In fact, here ordinary smoothing \cite{A} gives both uniform and $W^{2,2}$ strong approximation by an $\R^\ell$-valued smooth Sobolev function whose image lies in a small neighborhood of $N$; then composing this with the nearest-point projection to $N$ gives the desired smooth strong approximation with image in $N$. It was first observed in \cite{SU} that for $n=2$, a $W^{1,2}$ map (which may fail to have a continuous representative) admits strong $W^{1,2}$ approximation by smooth $W^{1,2}$ maps into $N$.  However for $n=3$, \cite{SU} also showed that the specific singular Sobolev map $x/|x|\in W^{1,2}(\B^3,\sS^2)$ is {\it not} strongly approximable in $W^{1,2}$ by smooth maps from $\B^3$ to $\sS^2$.  For first order Sobolev maps, the general problem of strong $W^{1,p}$ approximability was treated by F. Bethuel in \cite{Be2}, which (with \cite{BZ}) shows that: 
$$
W^{1,p}(\B^m,N)\ \mathrm{is\ the\ sequential\ strong}\ W^{1,p}\ \mathrm{closure\ of}\ {\cal C}^\infty(\B^m,N)\ \ \iff\ \  \Pi_{[p]}(N)=0\ .
$$
\noindent Here $[p]$ is the greatest integer less than or equal to $p$.  F. Hang and F. H. Lin, in \cite{HaL1} and \cite{HaL2}, updated these results with some new proofs and corrections, which account for the role played by the topology of the domain in approximability questions. See \cite{HaL2}, Th.1.3 for the precise conditions on the domain. There are many other interesting works on strong approximability of first order Sobolev maps by smooth maps, e.g. \cite{Be1}, \cite{BCL}, \cite{Hj},  \cite{BCDH}, \cite{BBC}, \cite{BZ}.
Generalization of the strong approximability results of \cite{Be2}, \cite{HaL1}, and \cite{HaL2} to higher order Sobolev mappings has been treated  by P. Bousquet, A. Ponce, and J. Van Schaftigen in \cite{BPV} which (with \cite{BZ}) shows that 
$$
W^{k,p}(\B^m,N)\ \mathrm{is\ the\ sequential\ strong}\ W^{k,p}\ \mathrm{closure\ of}\ {\cal C}^\infty(\B^m,N)\ \ \iff\ \  \Pi_{[kp]}(N)=0\ .
$$
\subsection{Sequential Weak Approximation}
The space $W^{k,p}(\Omega,N)$ also inherits the {\it weak topology} from $W^{k,p}(\Omega,\R^\ell)$.  A Sobolev map in $W^{k,p}(\Omega,N)$ that is not $W^{k,p}$ strongly  approximable by smooth maps may be $W^{k,p}$ {\it weakly approximable} by a  sequence  of smooth maps. For example, the map $x/|x|\in W^{1,2}(\B^3,\sS^2)$  {\it is} weakly approximable in $W^{1,2}$ by  some sequence $u_i$ of smooth maps. The well-known construction of such a $u_i$ involves changing  $x/|x|$ in a thin cylindrical tunnel $U$ of width $1/i$ extending from the origin $(0,0,0)$ to a point on $\partial\B^3$. To prove the weak $W^{1,2}$ weak convergence of $u_i$ to $x/|x|$, the key point of the construction is to keep the energies $\int_{\B^3}|\nabla u_i|^2dx$ bounded independent of $i$. 

To find an example of a map $v\in W^{1,p}(\B^m,N)$ which does not have a {\it weakly} approximating sequence of smooth maps, we need both $p<m$ and $\Pi_{[p]}(N)\neq 0$. Then,
{\it if $p$ is not an integer}, we simply choose, as in \cite{Be2}, any map $v$ which fails to have {\it strong} smooth approximations.  Assuming for contradiction that this map $v$ did admit some weak approximation by smooth maps $v_i$, then, for every point $a\in\B^m$, Fubini's theorem and Sobolev embedding (because $p>[p]$), would give {\it strong} convergence of the restrictions $v_i|\sS_a$ to $v|\sS_a$ for almost every $[p]$ dimensional Euclidean sphere $\sS_a$ centered at $a$ in $\B^m$. Then by the smoothness of $v_i$ and by [W], the corresponding homotopy classes $\lseg v|\sS_a\rseg$ would all vanish. But  Bethuel showed in [Be1] that precisely this local vanishing homotopy condition on $[p]$ spheres would imply that $v$ {\it does} admit strong smooth approximation, a contradiction. 

For integer $p$ the following question is still open:
\vskip.2cm
 {\it For any any compact manifold $N$, any integers $k,\,m\geq 1$ and any integer $p\geq 2$,  is \underline{every} Sobolev map $v\in  W^{k,p}(\B^m,N)$ actually  $W^{k,p}$ weakly approximable by a sequence of smooth maps?}
\vskip.2cm

This sequential weak density of smooth maps has been verified in the following cases:
\vskip.2cm
(1) \cite{BBC}, \cite{ABL} :\ \ $W^{1,p}(\B^m,\sS^p)$.
\vskip.1cm
(2) \cite{Hj}:\ \  $W^{1,p}(\B^m,N)$ with $N$ being simply $p-1$ connected (i.e. $\Pi_j(N)=0$ for $0\leq j\leq p-1$).
\vskip.1cm
(3) \cite{Pa}:\ \ $W^{1,1}(M,N)$ with $M$ and $N$ being arbitrary smooth manifolds with $\p N=\emptyset$ (weak convergence has to be understood in a {\it biting} sense here) 
\vskip.1cm
(4) \cite{PR}:\ \  $W^{1,2}(\B^m,N)$. (See also \cite{Ha} concerning the role of the topology of $M$ in  $W^{1,2}(M,N)$.)  
\vskip.2cm
See also a presentation of these results in \cite{Ri}. Another case is the main result of the present paper:
\vskip.2cm
\noindent\textbf {Theorem} \ref{wd}. {\it Any map in $W^{2,2}(\B^5,\sS^3)$ may be approximated in the  $W^{2,2}$ weak topology by a sequence of smooth maps}.
\vskip.2cm
In \S \ref{w22} below, we will explain how we came to study maps from $\B^5$ to $\sS^3$ and to look for $W^{2,2}$ estimates. But first we review a few of the ideas that were developed to study sequential weak convergence of smooth maps.
The space $W^{1,2}(\B^3,\sS^2)$ was studied extensively in the late eighties and early nineties with many works, e.g. \cite{HL}, \cite{BCL}, \cite{BZ}, \cite{BBC}, \cite{GMS1}. The concrete results of these many works has led to some analogous results and many conjectures for more general $k,n,p$, and $N$.  To sequentially weakly approximate a map $v\in W^{1,2}(\B^3,\sS^2)$, one first finds a $W^{1,2}$ {\it strong} approximation from the family $\mathcal R_0(\B^3,\sS^2)$ of maps $u\in  W^{1,2}(\B^3,\sS^2)$ which are smooth away from some finite set Sing$\,u$. In particular, we may assume $\int_{\B^3}|\nabla u|^2dx\leq 2\int_{\B^3}|\nabla v|^2dx$. Here the topology of $u$ near a point $a\in$Sing$\,v$ is given by the integer $d(a)=$degree\,[$u|\partial\B_\varepsilon(a)]$, which is independent of a.e. small $\varepsilon$. Then to get the desired completely smooth weak approximate, it is necessary to essential cancel the singularities of $u$. One does this by finding a one-chain or ``connection" $\Gamma_u$ with $\partial \Gamma_u$ in $\B^3$ being $\sum_{a\in\mathrm{Sing}\,v}d(a)\lseg a\rseg$ and with the rest of $\partial \gamma_u$ lying in $\partial B^3$. Then, as with the argument for $x/|x|$, one constructs smooth maps $u_i$ by making changes in tunnels of radius $1/i$ centered along the connection. To keep the $|\nabla u_i|^2$ integrals bounded, one needs to find a bound for the total length of the connection $\Gamma_u$ that depends only on $v$, and is independent of the approximating $u$. Here one may find a suitable connection by using the coarea formula.  This gives a good level curve of $u$ which connects the singularities to each other and to $\partial \B^3$ and which has length bounded by $\int_{\B^3}|\nabla u|^2dx$, which  has the independent bound $2\int_{\B^3}|\nabla v|^2dx$.

The first part of this argument, the strong $W^{1,2}$ approximation of an arbitrary  Sobolev map $v\in W^{1,2}(\B^3,\sS2)$ by a map $u\in\mathcal{R}_0(\B^3,\sS^2)$ has been generalized in \cite{Be2} to all $W^{1,p}(\B^m,N)$ and recently in \cite{BPV} to all $W^{k,p}(\B^m,N)$. Here one gets strong approximation by maps in $\mathcal R_{m-[p]-1}(\B^m,N)$ (respectively, $\mathcal R_{m-[kp]-1}(\B^m,N)$ which are smooth with singularities lying in finitely many affine planes of dimension $m-[p]-1$ (respectively, $m-[kp]-1$).

However, the second part involving canceling the singularities of $u$ has proven very challenging for generalization. One roughly needs an $m-[p]$ (respectively, $m-[kp]-1$) dimensional connection which has mass bounded in terms of the energy of $u$ and which the connects the singularity. Even with the connection, one still has to construct the bounded energy, smooth approximate.

\subsection{Topological Singularity and Bubbling}

In this supercritical dimension $m>kp$, we see that studying sequential $W^{k,p}$ weak smooth approximation in $W^{k,p}(\B^m,N)$ leads to questions about the relationship between the  possible energy drop, 
$\int_{\B^m}|\nabla^ku|^pdx < \liminf_{i\to\infty}\int_{\B^m}|\nabla^ku_i|^pdx$\ ,  of a $W^{k,p}$ weakly convergent sequence  $u_i\in W^{k,p}(\B^m,N)\cap \mathcal{C}^\infty$  and the possible singularities of its weakly convergent limit $u\in W^{k,p}(\B^m,N)$.  

For $0\neq\alpha\in\Pi_{kp}(N)$, we say a point $a\in\B^m$  is a type $\alpha$ {\it topological singularity} of a $W^{k,p}$ map $u$  if there is an $kp+1$ dimensional affine  plane $P$ containing $a$ so  the restrictions of $u$ to a.e. small $kp$ sphere $P\cap\partial\B_\varepsilon(a)$ induce (i.e. in the sense of \cite{W})  the homotopy class $\alpha$. Following the $W^{k,p}$ strong density of the partially smooth maps $\mathcal R_{m-kp-1}(\B^m,N)$ in $W^{k,p}(\B^m,N)$, one expects the  topological singularities, with their types as coefficients, to form a chain $S_u$ having dimension $m-kp-1$ and having coefficients in the group $\Pi_{kp}(N)$.  Recall the criterion of [Be1] that the vanishing of this ``$u$ topological singularity" chain( that is the vanishing of such homotopy classes for a.e. such restrictions at every $a\in\B^m$) is equivalent the $W^{k,p}$ strong approximability of $u$ by smooth maps.  

Also for $0\neq\alpha\in\Pi_{kp}(N)$,  the restrictions of $u_i$ to generic affine $kp$ planes can, as $i\to\infty$ have $|\nabla^k|^p$ energy concentration at an isolated point $b$ with an associated topological change corresponding to a type beta ``bubble''.  Putting such points together with their bubble types as coefficients should give a  

\noindent``$u_i$ bubbled" chain $B_{u_i}$ that has dimension $m-kp$, that has coefficient group $\Pi_{kp}(N)$, and that is carried by the $|\nabla^k(\,\cdot\,)|^p$ energy concentration set of the sequence.   

Using these vague definitions, one has the vague general conjecture:
\vskip.2cm
\noindent{\it\ \  Relative to  $\p\B^m$, the boundary of the $u_i$  bubbled chain $B_{u_i}$ equals the $u$ topological singularity chain $S_u$.} 
\vskip.2cm
\noindent The vagueness here concerns the precise definition of chain and boundary operation, and how one precisely obtains the bubbled chain $B_{u_i}$ from the sequence $u_i$ and the topological singular chain $S_u$ from $u$. From the cases we know, it is clear there is no single answer; it depends on the Sobolev space $W^{k,p}(\B^m,N)$, in particular the group $\Pi_{kp}(N)$.

In the special case $W^{1,2}(\B^3,\sS^2)$,  the relevant homotopy group is $\Pi_2(\sS^2)\simeq{\mathbb Z}$,  and \cite{BBC} and \cite{GMS1} show that this $1$ chain is precisely an integer-multiplicity 1 dimensional rectifiable current of finite mass (but possibly infinite boundary mass).  The special case was essentially generalized to $k=1$, $N=\sS^p$ in \cite{GMS2} and \cite{ABO}. Here,  the bubbled chain, now of dimension $m-p$, is again a rectifiable current. 

The paper [HR1] treated $W^{1,3}(\B^3,\sS^2)$. The relevant homotopy group is $\Pi_3(\sS^2)$, which is again isomorphic to ${\mathbb Z}$. Here the bubbled $1$ chain was shown to be possibly of infinite mass, and the notion of a ``scan'' was invented to describe precisely compactness and boundary properties. The paper [HR2] has able to handle bubbling in weak limits of smooth maps that corresponds to {\it any} nonzero homotopy class in the infinite nontorsion part of $\Pi_p(N)$. In this situation, the homotopy class of a map $w$ on the  sphere $\sS^{m-1}$ can again be described using a differential $m-1$ form $\Phi_w$ on $\sS^{m-1}$. The form is derived by a special algebro-combinatoric construction (depending on the rational homotopy class) involving a family of $w$ pullbacks of forms on $N$ and their ``$d^{-1}$ integrals". For example, in case $w:\sS^3\to\sS^2$, $\Phi_w = w^\#\omega_{\sS^2}\wedge d^*\Delta^{-1}w^\#\omega_{\sS^2}$. In general, this representation by a finite family of differential forms allows useful energy estimates involving certain Gauss integrals.  For a weakly convergent sequence of smooth maps, the bubbled chain, which cannot usually be represented as a finite mass current, can be understood precisely as a rectifiable scan whose boundary is given by the topological singularities of the limit Sobolev map. Though we have a somewhat satisfactory description of bubbling and topological singularity in all nontorsion cases, the question of sequential weak density in these cases are still not resolved, even for the case $W^{1,3}(\B^4,\sS^2)$.

Unfortunately representations of a homotopy class by differential forms are not available for {\it torsion} classes. In particular, if the relevant homotopy group of $N$ is completely torsion, then one requires other techniques to get energy estimates needed for questions about weak limits of smooth maps. The first such case is $W^{1,2}(\B^3,\R P^2)$, and was treated in [PR] . Here $\Pi_2(\R\pP^2)\simeq\Z_2$, and, a main result, is that {\it smooth maps are $W^{1,2}$ sequentially weakly dense} because $p=2$.

\subsection{$W^{2,2}(\B^5,\sS^3)$}\label{w22}

The present paper started with the modest goal of understanding analytic estimates for maps of

 \noindent$w:\sS^4\to \sS^3$ so as to understand weak convergence and sequential weak density for smooth maps from $\B^5$ to $\sS^3$.  Here, the appropriate homotopy group is $\Pi_4(\sS^3)$, which is isomorphic to $\Z_2$. But geometric descriptions of this homotopy class of $v$ are not very simple. As discussed in Sections \ref{pullback} and \ref{twisting} below, they involve considering, for a smooth approximation $u$ of $v$, the total twisting of a $u$-pullback normal framing upon circulation around a generic fiber $u^{-1}\{y\}$. The twisting of the normal  frame leads to an element of $\Pi_1\left(\sS\mathbb{O}(3)\right)\simeq \Z_2$.  To analytically compute such a twisting involves integration of a derivative  of a pull-back framing, hence a {\it second} derivative of the original map. So it is natural for this homotopy group to try to look for estimates in terms of the Hessian energy. 

A representative of the single nonzero element in  $\Pi_4(\sS^4,\sS^3)$ is the suspension of the Hopf map, $\sS\HH:\sS^5\to\sS^3\ $, described explicitly in the next section. In Section \ref{Prel} below, we slightly adapt \cite{BPV}, Th.5 by defining the subfamily
\begin{equation}\label{calR}
\mathcal{R} = \{u\in\mathcal{R}_0(\B^5,\sS^3)\ :\  u\,\equiv\, \sS{\HH}\lf(\frac{x-a}{|x-a|}\rg)\, \mathrm{on\ }\B_{\delta_0}(a)\setminus\{a\}\ \mathrm{for\,all\ } a\in\mathrm{Sing\,}u\  \mathrm{and\,some}\ \delta_0 > 0\}
\end{equation}
and then proving:
\vskip.3cm
\noindent{\bf Lemma \ref{Rdensity}}  {\it The family $\mathcal{R}$ is $W^{2,2}$ strongly dense in} $W^{2,2}(\B^5,\sS^3)$.
\medskip
\subsection{Connection Length} Given a finite subset $A$ of $\B^5$, one may define  $\Z_2$ {\it connection for $A$} (relative to $\partial\B^5$) as a finite disjoint union $\Gamma$ of finite length arcs embedded in $\overline{\B^5}$ whose union of endpoints is precisely $\ A\cup(\Gamma\cap\partial\B^5)\ $. Thus, each point of $A$ is joined by a unique arc in $\Gamma$ to either another point of $A$ or to a point of $\partial\B^5$. 

It will simplify some constructions to use a {\it minimal $\Z_2$ connection} for $A$, that is, one having least length. It is not difficult to verify the existence and structure of a  minimal $\Z_2$ connection for $A$. It simply consists of the disjoint union of finitely many closed intervals in $\B^5$ and finitely many radially pointing intervals having one endpoint in $\partial\B^5$. In \S \ref{sc}, we show how individual maps in $\mathcal R$ can be weakly approximated by smooth maps by proving: 

\smallskip
\noindent{\bf Theorem \ref{ML} (Singularity Cancellation)} {\it If $u\in\mathcal{R}$, $\Gamma$ is a minimal $\Z_2$ connection for Sing$\,v$, and $\varepsilon>0$, then there exists a smooth $u_\varepsilon\in W^{2,2}(\B^5,\sS^3)\cap\mathcal{C}^\infty$ so that $u_\varepsilon(x)=u(x)$ whenever dist$(x,\Gamma)>\varepsilon$ and
$$
\int_{\B^5}|\nabla^2u_\varepsilon|^2\,dx\ \leq\ \varepsilon\ +\ \int_{\B^5}|\nabla^2u|^2\,dx\ +\ c_{\sS\HH}\HF^1(\Gamma) 
$$
where} $c_{\sS\HH}\ =\ \int_{\sS^4}\lf|\nabla_{tan}^2(\sS{\HH})\rg|^2\,d\HF^4\ <\ \infty$.
\vskip.1cm
\noindent See Remark \ref{csSHH} concerning this constant.

\vskip.2cm

The core of our work, however, involves proving:

\vskip.2cm
\noindent{\bf Theorem \ref{ME} (Length Bound)} {\it For any $u\in\mathcal{R}$, Sing$\,u$ has a $\Z_2$ connection $\Gamma$ satisfying
$$
\HF^1(\Gamma)\ \leq\  c\int_{\B^5}|\nabla^2u|^2\,dx\ ,
$$ 
for some absolute constant $c$}.  
\vskip.2cm
Combining this length bound with Lemma \ref{Rdensity} and Theorem \ref{ML}, we readily establish, in Section \ref{sec3}, that any  Sobolev map in $W^{2,2}(\B^5,\sS^3)$ has a $W^{2,2}$
weak approximation by a sequence of smooth maps. 

We prove the length bound in Section \ref{length} by finding a suitable connection $\Gamma$ through three applications of the coarea formula. For a regular value $p\in\sS^3$ for $u$, the fiber $\Sigma=u^{-1}\{p\}$ is a smooth surface with cone point singularities at Sing$\,u$. By the coarea formula, we may choose this $p$ so that 
\begin{equation}\label{surfaceest}
\int_{u^{-1}\{p\}}\frac{|\nabla u|^4+|\nabla^2u|^2}{J_3u}\,d\HF^2\ .
\end{equation}
Then we need to choose connectiong curves on $\Sigma$. To do this we choose an orthonormal frame $\tilde\tau_1,\,\tilde\tau_2,\,\tilde\tau_3$ of the normal bundle of the surface $\Sigma=u^{-1}\{p\}$  by ortho-normalizing the $v$ pull-backs of a basis of Tan$\,(\sS^3,p)$. Inequality (\ref{surfaceest}) gives that 
$$
\int_\Sigma|\nabla\tilde\tau_j|\,d\HF^2\ \leq\ c\int_{\B^5}|\nabla^2u|^2\,dx\ . 
$$
We show how a a.e. oriented 2 plane in $\R^5$ determines at every point $x\in\Sigma$, with a finite exceptional set $b_1,\dots,b_j$, an orthogonal basis of Nor$(\Sigma,x)$, thought of as a reference normal framing. There is a unique $\gamma(x)\in\sS\oO(3)\simeq\R\pP^3$ and one gets some curves on $\Sigma$, with total length bounded by a multiple of $ c\int_{\B^5}|\nabla^2u|^2dx$ by choosing $\gamma^{-1}(E)$ where $E$ is a suitable great $\R\pP^2\subset\sS\oO(3)$. The curves starting at the some $a_i$ may end in either another $a_k$ or in $\p\B^5$ or (unfortunately) in a point $b_\ell$ where the reference framing degenerates. More argument, including another use of the coarea formula is required in sections \ref{3.8}, \ref{estimateB} to find additional curves of controlled length connecting $b_\ell$ to another $b_m$ or to $\p\B^5$. Putting all these curves together gives a $\Z_2$ connection for Sing$\,u$ satisfying the desired length bound.

\section{Preliminaries}\label{Prel}

We will let $c$ denote an absolute constant whose value may change
from statement to statement and which is usually  easily
estimable.  Here for $0\leq k \leq m$ and various $k$ dimensional subsets $A$ of $\R^m$, 
$$\int_Af\,d\HF^k\ =\ \int_Af(y)\,d\HF^ky
$$ 
will denote integration with respect to $k$ dimensional Hausdorff measure.  However, in top dimension where $\HF^m$ coincides with Lebesgue measure, we will use the use the standard notations $\int_A f\,dx\ =\ \int_A f(x)\,dx$.
\medskip
\begin{Lm}\label{Sobnorm} For each positive integer $m$, there is a positive constant $c_m$ so that
$$
\|v\|^2_{W^{2,2}(\B^m,N)}\ \leq\ c_m\left[\ (\mathrm{diam}\,N)^2\ +\ \int_{\B^m}|\nabla^2v|^2\,dx\ \right]
$$
for any compact Riemannian submanifold $N$ of $\R^\ell$ and $v\in W^{2,2}(\B^m,N)$.
\end{Lm}
{\it Proof.} Here $\|v\|^2_{W^{2,2}(\B^m,N)}\ =\ \int_{\B^m}\left(\,|v|^2+|\nabla v|^2+|\nabla^2v|^2\right)\,dx$. We clearly have the estimate
$$
\int_{\B^m}|v|^2dx\leq\ \HF^m(\B^m)(\mathrm{diam}\,N)^2\ .
$$
Moreover, by the Poincar\'e inequality,
$$
\begin{array}{lcl}
\ds\int_{\B^m}|\nabla v|^2\,dx\ &=&\ds\ \sum_{i=1}^m\int_{\B^m}\lf|\frac{\p v}{\p x_i}\rg|^2\,dx\leq\ \sum_{i=1}^m2\int_{\B^m}\lf|\frac{\p v}{\p x_i}-\lf(\frac{\p v}{\p x_i}\rg)_{avg}\rg|^2\,dx\ +\ 2\int_{\B^m}\lf|\lf(\frac{\p v}{\p x_i}\rg)_{avg}\rg|^2\,dx\\[5mm]
&\leq&\ds\ 2m\mathbf{C}_{\B^m}\int_{\B^m}|\nabla^2v|^2\,dx\ +\ \sum_{i=1}^m2\HF^m(\B^m)\lf|\lf(\frac{\p v}{\p x_i}\rg)_{avg}\rg|^2\ .
\end{array}
$$
It only remains to bound $\lf(\frac{\p v}{\p x_i}\rg)_{avg}$. We will do the case $i=1$, the cases $i\geq 2$ being similar. By Fubini's theorem and the absolutely continuity of $v$ on a.e. line in the $(1,0,\dots,0)$ direction,
$$
\begin{array}{l}
\ds\HF^m(\B^m)\lf|\lf(\frac{\p v}{\p x_1}\rg)_{avg}\rg|\ =\ \lf|\int_{\B^m}\frac{\p v}{\p x_1}\,dx\rg| \ =\ \lf|\int_{\B^{m-1}}\int_{-\sqrt{1-|y|^2}}^{\sqrt{1-|y|^2}}\frac{\p v}{\p x_1}(t,y_1,\dots,y_{m-1} )\,dt\,dy\rg| \\[5mm]
\ds\leq\ \int_{\B^{m-1}}\lf|v(\sqrt{1-|y|^2},y_1,\dots,y_{m-1} )-v(-\sqrt{1-|y|^2},y_1,\dots,y_{m-1})\rg| \,dy\ \leq\ \HF^{m-1}(\B^{m-1})\,\mathrm{diam}\,N\ .
\end{array}
$$
\cqfd
\vskip.2cm
\noindent\textbf{A Formula for the Suspension of the Hopf Map}
\vskip.2cm
Let $\HH:\sS^3\to \sS^2$ denote the standard Hopf map \cite{HR} :
$$
\HH(x_1,x_2,x_3,x_4)\ =\ \lf(\,2x_1x_2+2x_3x_4\,,\,2x_1x_4-2x_2x_3\,,\,x_1^2+x_3^2-x_2^2-x_4^2\, \rg) 
$$
 and ${\sS}{\HH}:\sS^4\to \sS^3$ be its suspension:
$$
{\sS}{\HH}(x_0,x_1,\cdots,x_4)\ =\
\lf(x_0\ ,\ \sqrt{1-x_0^2}\,\cdot\,\HH\lf(\frac{x_1}{\sqrt{x_1^2+\cdots +x_4^2}}\, ,\, 
\dots\, ,\, \frac{x_4}{\sqrt{x_1^2+\cdots+x_4^2}}\rg)\,\rg)\ .
$$
The latter map generates the nonzero
element of $\Pi_4(\sS^3)\simeq \Z_2$. Also, its homogeneous degree
$0$ extension
$$
\sS{\HH}(x/|x|)\ \in\ W^{2,2}(\B^5,\sS^3)\ .
$$
In particular, $\ \int_{\B^5}\lf|\nabla^2\lf(\sS{\HH}(x/|x|)\rg)\rg|^2dx\ =  c_{\sS\HH}\ $ where
\be\label{E}
c_{\sS\HH}\ =\ \int_{\sS^4}\lf|\nabla_{tan}^2(\sS{\HH})\rg|^2\,d\HF^4\ <\ \infty\ .
\ee
While the explicit formula for a suspension of the Hopf map is handy for simplifying proofs, the constant $c_{\sS\HH}$, which occurs in the conclusion of Theorem \ref{ML} can, by Remark \ref{csSHH}, be replaced by a more natural constant. 
\section{A Strongly Dense Family with Isolated Singularities}
Let ${\mathcal R}$ denote the class of $W^{2,2}(\B^5,\sS^3)$ maps that are
smooth except for finitely many suspension Hopf singularities.
That is,
$$
u\in{\mathcal R}\ \ \iff
$$
\be\label{delta0cond}
u\in\mathcal C^\infty(\B^5\setminus\{a_1,\dots,a_m\},\sS^3)\quad\mbox{ and }\quad  u(x)\ =\
\sS{\HH}\lf(\frac{x-a_i}{|x-a_i|}\rg)\quad\mbox{on}\quad \B_{\delta_0}(a_i)\setminus\{a_i\}
\ee
for some finite subset $\{a_1,\dots,a_m\}$ of $\B^5$ and some positive
$\delta_0<\min_i\{1-|a_i|, \min_{j\neq i}|a_i-a_j|/2\}\,\}$. 

\subsection{Strong Approximation by Maps in $\mathcal R$}

\begin{Lm}\label{Rdensity}  ${\mathcal R}$ is $W^{2,2}$ strongly dense in $W^{2,2}(\B^5,\sS^3)$.
\end{Lm}
{\it Proof.} Theorem 5 of \cite{BPV} gives the $W^{2,2}$ strongly density of the family
$\ \mathcal{R}_0^{2,2}(\B^5,\sS^3)\ $ of maps 

\noindent$v\in W^{2,2}(\B^5,\sS^3)$ which are smooth except for a finite singular set $\{a_1,\dots,a_m\}$ and which satisfy 
$$
\limsup_{x\to a_i}\lf(\ |x-a_i||\nabla v(x)|\ +\  |x-a_i|^2|\nabla^2v(x)|\ \rg)\ <\ \infty\ ,
$$
for $i=1,\dots,m$. Thus, it suffices to show: 

\vskip.4cm
{\it For each $\ v\in \mathcal{R}_0^{2,2}(\B^5,\sS^3)\ $ and $\varepsilon>0$, there is a map $\ u\in\mathcal R\ $ so that $\ \|u-v\|^2_{W^{2,2}} < \varepsilon\ $}.
\vskip.4cm

 Assuming Sing$\,v=\{a_1,\dots,a_m\}$, we will obtain $u$ by modifying $v$ near each point $a_i$. First fix a positive $\eta<\frac 12\min\lf\{\min_{i\neq j}|a_i-a_j|,\min_k(1-|a_k|)\rg\}$ so that 
\begin{equation}\label{Lbound}
L\ =\ \max_i\sup_{0<|x-a_i|<\eta}\lf(|x-a_i||\nabla v(x)|\ +\  |x-a_i|^2|\nabla^2v(x)|\rg)\ <\ \infty
\end{equation}
We will proceed in two stages: First we will find a positive $\delta<\eta$ depending only on 
$L$ and then define a map $w\in W^{2,2}(\B^5,\sS^3)$ so that 
\begin{equation}\label{w=v}
w\equiv v\ \mathrm{on}\ \B^5\setminus\cup_{i=1}^m\B_{\delta}(a_i)\ ,
\end{equation}
and, on each ball $\B_{\delta/2}(a_i)$, $w$ is degree-zero homogeneous about $a_i$, i.e.
$$
w(x)\ =\ w\lf(a_i + \frac{x-a_i}{|x-a_i|}\rg)\quad\mathrm{for}\quad0<|x-a_i|<\frac 12\delta\ .
$$
Second, we find a positive $\delta_0 <<\frac 12\delta$, depending on $w|\cup_{i=1}^m\p\B_{\delta/2}(a_i)$, and a map $u\in\mathcal R$ with $u\equiv w$ on 

\noindent$\B^5\setminus\cup_{i=1}^m\B_{\delta_0}(a_i)$ and, on each ball $\B_{\delta_0}(a_i)$,
$u\equiv\sS\HH\lf(\frac{x-a_i}{|x-a_i|} \rg)$.

For the first step, we first fix a smooth monotone increasing $\lambda:[0,\infty)\to [\frac 12,\infty)$ so that
$$
\lambda(t)\ =\ \cases{1/2 & for\quad $0\leq t\leq\frac 12$\cr
t & for\quad $t\geq 1$\ .\cr}
$$

Consider the unscaled situation of a map $V\in\mathcal{C}^\infty(\overline{\B^5}\setminus\B_{\frac 12}(0),\sS^3)$ with 

\noindent$|\nabla V|+|\nabla^2V|\leq L$.  Then we define the reparameterized map 
$$
W(x)=V\lf(\lambda(|x|)x\rg)\ ,
$$
and see that, with respect to the radial variable 
$\rho=|x|$, 
$$
\frac{\p W}{\p \rho}\ \equiv\ \frac{\p V}{\p \rho}\quad \mathrm{on}\quad \p\B^5\ ,\quad\quad
\frac{\p W}{\p \rho}\ \equiv\ 0\quad \mathrm{on}\quad \overline{\B_{\frac 12(0)}}\setminus\{0\}\ ,
$$
Moreover, using explicit pointwise bounds for $|\lambda'|$ and $|\lambda''|$,  we readily find an explicit  constant $C$ so that 
$$
|\nabla W(x)|+|\nabla^2W(x)|\ \leq\ C\,L\quad \mathrm{for}\quad \frac 12\leq|x|\leq 1\ ,
$$
$$
|x||\nabla W(x)|+|x|^2|\nabla^2W(x)|\ \leq\  C\,L\quad \mathrm{for}\quad  0<|x|< \frac 12\ \ .
$$
\vskip.1cm  
Now we return to the original scale by defining $w$ to satisfy (\ref{w=v}) and to have, for each point $x\in\B_{\delta}(a_i)$,
$$
w(x)\ =\ W\lf(\frac{x-a_i}{\delta}\rg)\quad\mathrm{where}\quad
V(x)\ =\ v(a_i+\delta x)\ .
$$
Then $w$ belongs to $W^{2,2}(\B^5,\sS^3)$ and satisfies the estimate
$$
\max_i\sup_{0<|x-a_i|<\delta}\lf(|x-a_i||\nabla w(x)|\ +\  |x-a_i|^2|\nabla^2w(x)|\rg)\ \leq\ CL\ ,
$$
hence,
$$
\begin{array}{lcl}
\ds\|w-v\|^2_{W^{2,2}} &=& \sum_{i=1}^m\int_{\B_{\delta}(a_i)}\lf(|w-v|^2+|\nabla w-\nabla v|^2+|\nabla^2w-\nabla^2v|^2\rg)dx\\[5mm] 
\ds&\leq& 2\sum_{i=1}^m\int_{\B_{\delta}(a_i)}\lf(|w|^2+|v|^2+|\nabla w|^2+\nabla v|^2+|\nabla^2w|^2+|\nabla^2v|^2\rg)dx\ \leq\ c(1+L)^2\delta\ .
\end{array}
$$
So we easily choose $\delta$ so that $\|w-v\|^2_{W^{2,2}}<\frac 14\varepsilon$.
\vskip.2cm
For the second step we note that any continuous homotopy between smooth maps of smooth manifolds can be made smooth; hence\ :  
$$
\mathit{A\ smooth\ map\ }\phi:\sS^4\to\sS^3\ \mathit{is\ smoothly\ homotopic}
\cases{\mathit{either\ to\ a\ constant} & \it{in case} $\lseg \phi\rseg = 0\in \Pi_4(\sS^4,\sS^3)$\cr
\mathit{or\ to\quad } \lseg\sS\HH\rseg & \it{in case} $\lseg \phi\rseg \neq 0\in \Pi_4(\sS^4,\sS^3))$.\cr} 
$$

For each $i=1,\dots,m$, we apply this to the map $\phi_i(x)=w(a_i+ \frac 12\delta x)$ to obtain a smooth homotopy $h_i:[0,1]\times\sS^4\to\sS^3$ which connects $\phi_i$ to $\sS\HH$. Reparameterizing the time variable near $0$ and $1$, we may assume 
$$
h_i(t,y)\ =\ \cases{ \phi_i(y) &for $t$ near\ $0$\cr
(\sS\HH)(y) &for $t$ near\ $1$\ .\cr}
$$
 By smoothness, 
$K\ =\ \sup_i\| h_i\|_{W^{2,2}}<\infty$, and we will, for some $\delta_0<<\frac 12\delta$, define the map $u$ by 
$$
u\equiv w\quad \mathrm{on}\quad \B^5\setminus\cup_{i=1}^m\B_{\delta_0}(a_i)\ ,
$$
$$
u(x)\ =\ h_i\lf(2-2\frac{|x-a_i|}{\delta_0}\ ,\ \frac{x-a_i}{\delta_0}\rg)\quad \mathrm{for}\quad x\in\B_{\delta_0}(a_i)\setminus \B_{\frac 12\delta_0}(a_i)\ ,
$$
$$
u(x)\ =\ \sS\HH\lf(\frac{x-a_i}{|x-a_i|}\rg)\quad \mathrm{for}\quad x\in\B_{\frac 12\delta_0}(a_i)\setminus\{0\}\ , 
$$
for $i=1,\dots,m$.  One readily checks that $u\in\mathcal{R}$. Moreover, as in Step 1, we find that 
$$
\|u-w\|^2_{W^{2,2}}\ = \sum_{i=1}^m\int_{\B_{\delta}(a_i)}\lf(|u-w|^2+|\nabla u-\nabla w|^2+|\nabla^2u-\nabla^2w|^2\rg)dx\ \leq\ c(1+K)^4\delta_0\ .
$$
So we easily choose $\delta_0$ small enough so that $\ \|u-w\|^2_{W^{2,2}}<\frac 14\varepsilon\ $ and obtain the desired estimate
$$
\|u-v\|^2_{W^{2,2}}\ \leq 2\|u-w\|^2_{W^{2,2}}\ +\ 2\|w-v\|^2_{W^{2,2}}\ <\ \varepsilon\ .  
$$
\cqfd

\subsection{Insertion of an $\sS\HH$ Bubble into a Map from $\B^4$ to $\sS^3$}\label{binsert}

Arguing as in the proof of Lemma \ref{Rdensity}, we first fix a monotone increasing smooth function $\mu$ on $[0,\infty)$ so that
$$
\mu(t)\ =\ \cases{0 & for\quad $0\leq t\leq\frac 12$\cr
t & for\quad $t\geq 1$\ ,\cr}
$$  
$|\mu'|\leq 3$,
and $|\mu''|\leq 16$.  We readily prove the following:

\begin{Lm}[Initial Reparameterization]\label{fsig} There an absolute constant $C$ so that, for any smooth $f :\B^4\to\sS^3$ and $0<\sigma<1$, the map
$$
f_\sigma :\B^4\to\sS^3\ ,\quad f_\sigma(y)\ =\ f\lf(\mu\lf(\sigma^{-1}|y|\rg)|y|^{-1}\sigma y\rg)\quad\mathrm{for}\ y\in\B^4\ ,
$$
coincides with $f$ on $\B^4\setminus\B^4_{\sigma}$, is identically equal to $f(0)$ on $\B^4_{\sigma/2}$, and satisfies
$$
\sup_{\B^4_\sigma}|\nabla f_\sigma|\ \leq\ C\sup_{\B^4_\sigma}|\nabla f_\sigma|\ ,\quad\quad\sup_{\B^4_\sigma}|\nabla^2f_\sigma|\ \leq\ C\sigma^{-1}\sup_{\B^4_\sigma}|\nabla^2f| \ .
$$
In particular,
\be\label{resest}
\int_{\B^4_\sigma}|\nabla^2f_\sigma |^2\,dx\ \leq\ C^2\HF^4(\B^4)\sup_{\B^4_\sigma}|\nabla^2f |^2\,\sigma^2\ .
\ee
\end{Lm}
\vskip.2cm
\noindent\textbf{Construction of an $\sS\HH$ Bubble}
\vskip.1cm
Recalling that $\sS\HH(0,1,0,0,0) = (0,0,0,1)$, we will first slightly modify $\sS\HH$ to be constant near $(0,1,0,0,1)$. Consider the spherical coordinate parameterization,
$$
\Upsilon : [0,\pi]\times\sS^3\to\sS^4\ ,\quad\Upsilon(\rho,\omega) = \lf( (\sin \rho)\omega_1,\cos \rho, (\sin \rho)\omega_2, (\sin \rho)\omega_3,(\sin\rho)\omega_4\rg)\ .
$$
Arguing again as in the proofs of Lemma \ref{Rdensity} and Lemma \ref{fsig},  we let  
$$
M_{\rho_0}(\rho,\omega)\ =\ \lf(\mu(\rho_0^{-1}\rho)\rho_0\rho\,  , \omega\rg)\quad\mathrm{for}\quad 0\leq\rho\leq\rho_0<<1\quad\mathrm{and}\quad\omega\in \sS^3\ ,
$$
and define $\Phi_{\rho_0}:\sS^4\to\sS^4$ by
$$
\Phi_{\rho_0}(y)\ =\ \cases{ 
\Upsilon\circ M_{\rho_0}\circ\Upsilon^{-1}\{y\}  &  for $y\in\Upsilon\lf( [0,\rho_0 ]\times\sS^4\rg)$\cr
$y$ & otherwise\ .\cr}
$$
Then $\Phi_{\rho_0}$ is surjective, and maps the entire spherical cap 
$$
\Omega_{\rho_0/2}\ =\ \Upsilon\lf( [0,\rho_0/2) \times\sS^4\rg)\ =\ \sS^4\cap\B^5_{2\sin(\rho_0/4)}\lf( (0,1,0,0,0)\rg)
$$
 to its center point $(0,1,0,0,0)$. 

We now consider the composition $\sS\HH\circ\Phi_{\rho_0}$   The homotopy class is unchanged
\be\label{homotopy}
\lseg\sS\HH\circ\Phi_{\rho_0}\rseg\ =\ \lseg\sS\HH\rseg\ \neq 0\ \in\ \Pi_4(\sS^3)\ .
\ee
Noting  that $\rho_0\frac{\p^2 M_{\rho_0}}{\p\rho^2}$ is bounded independent of $\rho_0$ and has support in $\Upsilon\lf([0,\rho_0]\times\sS^3\rg)$, we readily verify, as in (\ref{resest}), that
\be\label{modShopf}
\int_{\sS^4}|\nabla_{tan}^2(\sS\HH\circ\Phi_{\rho_0})|^2d\HF^4\
 =\ \int_{\sS^4\setminus\Omega_{\rho_0}}|\nabla_{tan}^2(\sS\HH )|^2d\HF^4 +\ \int_{\Omega_{\rho_0}}|\nabla_{tan}^2(\sS\HH\circ\Phi_{\rho_0})|^2d\HF^4
\to\quad  c_{\sS\HH}\ +\ 0\  ,
\ee
as $\rho_0\to 0$.

Since the stereographic projection 
$$
\Pi:\sS^4\setminus\{(0,1,0,0,0)\}\to\R^4\ ,\quad\Pi(x_0,x_1,x_2,x_3,x_4)\ =\ \lf(\frac{x_0}{1-x_1},\frac{x_2}{1-x_1},\frac{x_3}{1-x_1},\frac{x_4}{1-x_1}\rg)
$$
is conformal,  we  get the conformal diffeomorphism
$$
\Lambda_{\rho_0}\ :\ \overline{\B^4}\ \to\ \sS^4\setminus \Omega_{\rho_0/2}\ ,\quad\Lambda_{\rho_0}(y)\ =\ \Pi^{-1}\lf[\lf(\frac{\sin(\rho_0/2)}{1-\cos(\rho_0/2)}\rg)\,y\rg]\ .
$$
We see that 
$$
\sS\HH\circ\Phi_{\rho_0}\circ \Lambda_{\rho_0}\ :\ \overline{\B^4}\ \to\ \sS^3
$$
is a smooth surjection which sends $\p\B^4$ identically to $\sS\HH(0,1,0,0,0) = (0,0,0,1)$. 

We wish to have a similar map that has boundary values being a a possibly different constant $\xi\in\sS^3$. To get a suitable formula, it will be handy to recall that $\sS^3$, being identified with the unit quaternions,
$$
\{\,x_1+x_2\mathbf{i}+x_3\mathbf{j}+x_4\mathbf{k}\ :\ x_1^2+x_2^2+x_3^2+x_4^2=1\,\}
$$
has a product $ \star\ $ and, in particular, that 
$$
 (0,0,0,-1)\star(0,0,0,1)\ =\ (-{\bf k})\star{\bf k}\ = {\bf 1}\ =\ (1,0,0,0)\ .
$$
So now we define, for every  $\rho_0>0$ and  $\xi\in\sS^3$, the bubble
\be\label{bubbledef}
\mathcal{B}_{\rho_0,\xi}\ :\ \B^4\ \to\ \sS^3\ ,\quad\mathcal{B}_{\rho_0,\xi}\ \equiv\ \xi\star (-\mathbf{k})\star\sS\HH\circ\Phi_{\rho_0}\circ \Lambda_{\rho_0}\ ,
\ee
which is identically equal to $\xi$ on $\p\B^4$.  By the conformal invariance of the Hessian energy and (\ref{modShopf}),
\be\label{bblenergy}
\int_{\B^4}|\nabla^2\mathcal{B}_{\rho_0,\xi}|^2\,dy\ =\ \int_{\sS^4}|\nabla^2_{tan}(\sS\HH\circ\Phi_{\rho_0})|^2\,d\HF^4\ \to\ \int_{\sS^4}|\nabla_{tan}^2(\sS\HH)|^2\,d\HF^4\ =\ c_{\sS\HH}\quad\quad\mathrm{as}\ \rho_0\to 0\ .
\ee

Now, for a smooth map $f:\B^4\to\sS^3$, $0<\sigma<1$,  and $0<\rho_0 <<1$, we define
the ``bubbled'' map $f_{\sigma,\rho_0}:\B^4\to\sS^3$,
\be
\label{fwithbubble}
f_{\sigma,\rho_0}(y)\ =\ \cases{ f_{\sigma} (y) &  for $\sigma/2< |y|< 1$\ , \cr  \mathcal{B}_{\rho_0,f(0) }(2y/\sigma) &  for $|y|\leq\sigma/2$ .\cr }
\ee
\vskip.2cm
\noindent Note that $f_{\sigma,\rho_0}$ is continuous because both expressions equal $f(0)$ on $\p\B^4_{\sigma/2}$. Moreover, all the positive order derivatives of both expressions vanish here because the smooth map $f_\sigma$ is  constant on $\B^4_{\sigma/2}$ while the smooth reparameterization $\Phi_{\rho_0}$ is constant on $\Omega_{\rho_0/2}$. So the map $f_{\sigma,\rho_0}$ is, in fact, smooth.  Moreover, since the Hessian energy is invariant under change of scale,
\be\label{fsrint}
\begin{array}{lcl}
\ds \int_{\B^4}\lf|\nabla^2f_{\sigma,\rho_0}\rg|^2dy &=& \int_{\B^4\setminus\B^4_{\sigma/2}}\lf|\nabla^2f_{\sigma}\rg|^2dy\ +\  \int_{\B^4_{\sigma/2}}\lf|\nabla^2\lf[{\mathcal{B}}_{\rho_0,f(0) }\lf(\frac 2{\sigma}(\cdot )\rg) \rg]\rg|^2dy\\[5mm] 
\ds  &=&\ \int_{\B^4\setminus\B^4_\sigma}\lf|\nabla^2f\rg|^2\,dy\ +\int_{\B^4_\sigma\setminus\B^4_{\sigma/2}}\lf|\nabla^2f_{\sigma}\rg|^2\,dy\ +\  \int_{\B^4}\lf|\nabla^2{\mathcal{B}}_{\rho_0,f(0) }\rg|^2\,dy\\[5mm] 
\ds&\to&\ \int_{\B^4}\lf|\nabla^2f\rg|^2\,dy\quad +\quad  0\quad +\quad \int_{\B^4}\lf|\nabla^2\mathcal{B}_{\rho_0,f(0) }\rg|^2\,dy\quad\quad \mathrm{as}\quad \sigma\to 0\ .
\end{array}
\ee
Note that the middle equation also gives, with (\ref{resest}) the bound
\be\label{fbd}
 \int_{\B^4}\lf|\nabla^2f_{\sigma,\rho_0}\rg|^2dy\ \leq\ (1+C^2)\sup_{\B^4}|\nabla^2 f|^2\ + \int_{\B^4}\lf|\nabla^2\mathcal{B}_{\rho_0,f(0) }\rg|^2\,dy ,
\ee
independent of $\sigma$.

\subsection{Singularity Cancellation of a Map in $\mathcal R$}\label{sc} 
\vskip.2cm
\begin{Th}\label{ML} If $u\in\mathcal{R}$, $\Gamma$ is a minimal $\Z_2$ connection for Sing$\,v$, and $\varepsilon>0$, then there exists a smooth $u_\varepsilon\in W^{2,2}(\B^5,\sS^3)\cap\mathcal{C}^\infty$ so that $u_\varepsilon(x)=u(x)$ whenever dist$(x,\Gamma)>\varepsilon$ and
$$
\int_{\B^5}|\nabla^2u_\varepsilon|^2\,dx\ \leq\ \varepsilon\ +\ \int_{\B^5}|\nabla^2u|^2\,dx\ +\ c_{\sS\HH}\HF^1(\Gamma)\ . 
$$
\end{Th}
{\it Proof.}   Note that, by slightly rescaling near $\p\B^5$, we may assume that $u$ extends smoothly to a neighborhood of $\overline{\B^5}$.

First, using (\ref{bblenergy}), we fix  a positive $\rho_0$ to be small enough so that, 
\be\label{rho0}
\lf| c_{\sS\HH}\ -\ \int_{\B^4}|\nabla^2\mathcal{B}_{\rho_0,\xi}|^2\,dy\ \rg|\ <\ \frac\varepsilon{4(1+\HF^1(\Gamma))}\ 
\ee
for all $\xi\in\sS^3$.  Also we recall that a minimal $\Z_2$ connection for $u$ is the union of a finite family $\mathcal{I}$ of disjoint closed intervals $I=[a_I,b_I]$ where
\be\label{aIbI}
a_I\ \in\ \mathrm{Sing}\,u\quad\mathrm{and}\quad b_I\ \in\ \cases{either\ &$\mathrm{Sing}\,u$\cr
                                   or\ &$\p\B^5$ with $I\perp\p\B^5$\ .\cr} 
\ee
\vskip.3cm
Second, we fix a positive $\delta_0$ so that:
\begin{itemize}

\item[(1)] $\delta_0<\frac 12\min_{a\in\mathrm{Sing}\,u}\{1-|a|, \min_{a \neq \tilde a\in\mathrm{Sing}\,u}|a-\tilde a|\,\}$

\item[(2)] $\delta_0\ < \frac 12\min\{ |x-\tilde x|\ :\ x\in I,\ \tilde x\in\tilde I,\ I\neq\tilde I\in\mathcal{I}\ \}$. 

\item[(3)] $u\in\mathcal{C}^\infty(\B^5_{1+\delta_0}\setminus\mathrm{Sing}\,u,\sS^3)$.  

\item[(4)]$\delta_0<(1+c_{\sS\HH})^{-1}\lf(1+\mathrm{card}(\mathrm{Sing}\,u)\rg)^{-1}\varepsilon/5$
\end{itemize}
Our main step in constructing $u_\varepsilon$ will be to use, for each $I\in\mathcal I$,  the  bubble insertion of \S \ref{binsert} in each cross-section of a pinched cylindrical region $V_I$  of radius $\delta_0/9$. Near the singular endpoints of $I$, $V_I$ is pinched to be a round cone with opening angle $2\arctan(\frac 19)$.  

To describe the explicit construction, we need some notation.  With $I=[a_I,b_I]\in\mathcal{I}$ as above in (\ref{aIbI}), let 
$$
|I|\ =\ |b_I-a_I|,\quad\mathbf{e}_I\ =\ (b_I-a_I)/|I|\ \in\sS^4 ,\quad\quad\Pi_I:\R^5\to\R,\quad\Pi_I(x)=x\cdot\mathbf{e}_I-a_I\cdot\mathbf{e}_I\ ,
$$
and $ B_I(t,r) $ be the open  ball in the $4$ dimensional affine plane $\Pi_I^{-1}\{t\}$ with center $a_I+t\mathbf{e}_I$ and radius $r$.  We now define $\ V_I\ $ to be the pinched cylindrical region
$$
V_I\ =\ \bigcup_{0<t<|I|} B_I\lf(t, r_I(t)\rg)\ ,
$$
by using a fixed smooth function $\nu: [0,\infty)\to [0,1]$ with 
$$
\nu(t)\ =\ \cases{ t & for $0\leq t\leq \frac 12$\cr
1 & for $1\leq t$\ ,\cr}
$$ 
to define the smooth radius function 
$$r_I(t)\ =\ \cases{ (\delta_0/9)t &for $0\leq t \leq \frac 12 \delta_0$\cr
(\delta_0/9)\nu(t/\delta_0) & for $\frac12 \delta_0\leq t\leq\delta_0$\cr
\delta_0/9 & for $\delta_0\leq t\leq |I|-\delta_0$\cr
(\delta_0/9)\nu((|I|-t)/\delta_0) & in case $b_I\in\mathrm{Sing}\,u$ and $|I|-\delta_0\leq t\leq |I|-\frac 12\delta_0$\cr
(\delta_0/9)\,(|I|-t) & in case $b_I\in\mathrm{Sing}\,u$ and  $|I|-\frac 12\delta_0\leq t\leq |I|$\cr
\delta_0/9 & in case $b_I\in\p\B^5$ and  $|I|-\delta_0\leq t\leq |I|$\ .\cr}
$$
Thus $\p V_I$ is a smooth hypersurface except for the conepoint(s) $a_I$ (and $b_I$ in case $b_I\in\mathrm{Sing}\,u$).

For convenience, we fix an orthonormal basis $\mathbf{e}^I_1,\dots,\mathbf{e}^I_4$ for the orthogonal complement of $\ \R\mathbf{e}_I$. 
For each  $t\in\R$, we will use the affine similarity 
$$
A_{I,t}:\R^4\to \R^5\ ,\quad A_{I,t}(y_1,\dots ,y_4)\ =\ a_I + t\mathbf{e}_I + r_I(t)\lf[ y_1\mathbf{e}^I_1+\cdots +y_4\mathbf{e}^I_4\rg]
$$ 
so that $\ A_{I,t}(\B^4)\ =\ B_I\lf(t,r_I(t)\rg)$. 

For $0<\sigma<1$, we recall (\ref{fwithbubble}) and  define the smooth reparameterized map
$$
v_\sigma(x)\ = \cases{u(x) & for $x\in\B^5\setminus\bigcup_{I\in\mathcal{I}}V_I\cup\{a_I\}\cup\{b_I\}$\cr
	\lf(u\circ A_{I,t}\rg)_{\sigma,\rho_0}\lf(A_{I,t}^{-1}(x)\rg) & for $x\in B_I\lf(t, r_I(t)\rg)$\ .\cr}
$$
Observe that $ v_\sigma$ actually coincides with $u(x)$ outside the ``$\sigma$ thin'' set $\ \bigcup_{I\in\mathcal{I}}V^\sigma_I\cup\{a_I\}\cup\{b_I\}\ $ where
$$
V^\sigma_I\ =\ \bigcup_{0<t<|I|} B_I\lf(t,\sigma r_I(t)\rg)\ .
$$

The explicit formulas given above and in the earlier parts of \S \ref{sc} show the qualitative smoothness of $v_\sigma$ on $\B^5\setminus\mathrm{Sing}\,u$. Our next goal is to verify that 
\be\label{wsigest}
\limsup_{\sigma\to 0}\sum_{I\in\mathcal{I}}\int_{V_I^\sigma}\lf|\nabla^2 v_\sigma\rg|^2\,dx\quad <\quad c_{\sS\HH}\HF(\Gamma)\ +\ \frac \varepsilon 2\ .
\ee

We define
$$
J_I\ =\ \lf\{t\ :\ \frac 12\delta_0\ \leq\ t\ \leq\ \cases{ |I|-\frac 12\delta_0 & in case\ $b_I\in\mathrm{Sing}\,u$\cr
|I|\ & in case $b_I\in\p\B^5$\ ,\cr}\ \rg\}
$$
and note that on $J_I$, the scaling factor $r_I(t)$ satisfies $\frac 1{18}\delta\leq r_I(t)\leq \frac 19\delta_0$ while $|r_I'(t)|$ and $|r_I''(t)|$ are bounded.  We will use the corresponding truncated sets
$$
W_I\equiv V_I\cap\Pi_I^{-1}(J_I)\ =\ \bigcup_{t\in J_I} B_I\lf(t, r_I(t)\rg)\ ,\quad
W^\sigma_I\equiv V^\sigma_I\cap\Pi_I^{-1}(J_I)\ =\ \bigcup_{t\in J_I} B_I\lf(t, \sigma r_I(t)\rg)\ .
$$
We have the pointwise bound 
\be\label{supbd}
L\ =\ \sup_{W_I}\lf(|\nabla u|^2+|\nabla^2u|^2\rg)\ <\ \infty\ ;
\ee
hence, by (\ref{fbd}),
$$
\sup_{t\in J_I}\sup_{0<\sigma<1}\int_{\B^4}\lf|\nabla^2(u\circ A_{I,t})_{\sigma,\rho_0}\rg|^2\,dy\ <\ \infty\ .
$$

Note that the orthogonality of the five vectors $\mathbf{e}_I, \mathbf{e}^I_1,\cdots,\mathbf{e}_4^I$  lead to the decomposition of the squared  Hessian norm into pure second partial derivatives
$$
|\nabla^2( \cdot)|^2\quad =\quad \lf|\nabla_{\mathbf{e}_I,\mathbf{e}_I}( \cdot)\rg|^2\ +\ \lf|(\nabla_{\mathbf{e}^I_1,\mathbf{e}^I_1}(\cdot)\rg|^2\  +\ \cdots\ +\ \lf|(\nabla_{\mathbf{e}^I_4,\mathbf{e}^I_4}(\cdot)\rg|^2\ ,
$$
which we will abbreviate as $\lf|\nabla^2_{\mathbf{e}_I}( \cdot)\rg|^2\ +\ \lf|\nabla^2_{\mathbf{e}_I^{\perp}}( \cdot)\rg|^2\ $.

 It follows from Fubini's Theorem,  the conformal invariance of the $4$ dimensional Hessian energy, (\ref{fsrint}),   dominated convergence, and (\ref{rho0})  that 
\be\label{transverse}
\begin{array}{lcl}
\ds \int_{W^\sigma_I}\lf|\nabla^2_{\mathbf{e}_I^{\perp}}v_\sigma\rg|^2\,dx\ &=&
\int_{t\in J_I}\int_{B_I(t,r_I(t))}\lf|\nabla^2_{\mathbf{e}_I^{\perp}}v_\sigma\rg|^2\,d\HF^4\,dt
\\[5mm] 
\ds &=& 
\int_{t\in J_I}\int_{\B^4}\lf|\nabla^2\lf(u\circ A_{I,t}\rg)_{\sigma,\rho_0}\rg|^2\, dy\,dt
\\[5mm] 
\ds &\to& 
\int_{t\in J_I}\lf[\ 0\ +\int_{\B^4}\lf| \nabla^2\mathcal{B}_{\rho_0,u\lf(a_I+t\mathbf{e}_I\rg) }\rg|^2\,dy\rg]\,dt\ \leq\   c_{\sS\HH}|I|\ +\ \frac{\varepsilon |I|}{4(1+\HF^1(\Gamma ))}\ ,
\end{array}
\ee
as $\sigma\to 0$.  Thus 
\be\label{wsigperp}
\limsup_{\sigma\to 0}\sum_{I\in\mathcal{I}}\int_{W^{\sigma}_I}\lf|\nabla_{\mathbf{e}_I^{\perp}}^2v_{\sigma}\rg|^2\,dx\quad \leq\quad \sum_{I\in\mathcal{I}}\lf[c_{\sS\HH}|I| +\frac{\varepsilon |I|}{4(1+\HF^1(\Gamma ))}\rg]\quad \leq\quad c_{\sS\HH}\HF^1(\Gamma)\ +\ \varepsilon/4\ .
\ee

To get the full squared Hessian integral $\int_{W^\sigma_I}\lf|\nabla^2v_\sigma\rg|^2\,dx$, we also need to consider $\int_{W^\sigma_I}\lf|\nabla^2_{\mathbf{e}_I}v_\sigma\rg|^2\,dx$, which involves computing $\frac {\p^2}{\p t^2}$ of various terms.  To estimate the last integral, it again suffices by Fubini's theorem, Lemma \ref{fsig}, (\ref{fwithbubble}), and changing variables to consider
\be\label{2terms}
\int_{t\in J_I}\lf[\ \int_{\B^4_\sigma\setminus\B^4_{\sigma/2}}\lf| \frac {\p^2}{\p t^2}(u\circ A_{I,t})\rg|^2\,dy\,dt\ +\ \int_{\B^4_{\sigma/2}}\lf|\frac {\p^2}{\p t^2}\lf[\mathcal{B}_{\rho_0,u\lf(a_I+t\mathbf{e}_I\rg) }\lf(\frac 2\sigma(\cdot)\rg)\rg]\ \rg|^2\,dy\ \rg]\ \,dt\ .
\ee 
The chain rule and the bounds of $|r'|$ and $|r''|$ on $J_I$ give the pointwise bound
$$
\lf|\frac {\p^2}{\p t^2}(u\circ A_{I,t})\rg|\ =\ \lf|\frac {\p^2}{\p t^2}\lf[u\lf( a_I + t\mathbf{e}_I + r_I(t)y\rg) \rg] \rg|\ \leq\ cL\ ,
$$
and definition (\ref{bubbledef})  gives the bound
$$
\lf|\frac {\p^2}{\p t^2}\mathcal{B}_{\rho_0,u\lf(a_I+t\mathbf{e}_I\rg)}\rg|\ =\ \lf|\frac {\p^2}{\p t^2}\lf[u\lf(a_I+t\mathbf{e}_I\rg)\star (-\mathbf{k})\star\sS\HH\circ\Phi_{\rho_0}\circ \Lambda_{\rho_0}\ \rg]\rg|\ \leq\ L\ .
$$
Integrating implies that  (\ref{2terms}) is bounded by $\ (c+1)L|I|\HF^4(\B^4_{\sigma})\ $, and we deduce that 
\be\label{tintegral}
\lim_{\sigma\to 0}\int_{W^\sigma_I}\lf|\nabla^2_{\mathbf{e}_I}v_\sigma\rg|^2\,dx\ =\ 0\ .
\ee

Next we consider the conical end(s) $V_I^\sigma\setminus W_I^\sigma$. By our choice of $\delta_0$, $u$ is degree-$0$ homogeneous about $a_I$ on the region $\B^5_{\delta_0}(a_I)$.  It follows that all of the normalized bubbled functions $(u\circ A_{I,t})_{\sigma,\rho_0}$ coincide for $0<t\leq\sigma/2$. Thus, in the one conical end $\ V_I^\sigma\cap\Pi_I^{-1}(0,\delta_0/2]\ $, $v_\sigma$ is also degree-$0$ homogeneous about $a_I$. So we can easily estimate the Hessian integral there by using spherical coordinates about $a_I$. Note that radial projection of the $4$ dimensional Euclidean ball $\ V_I^\sigma\cap\Pi^{-1}\{\delta_0/2\}\ $ onto the small spherical cap $V_I^\sigma\cap\p\B^5_{\delta_0/2}$ is a smooth diffeomorphism with easily computed $\mathcal{C}^2$ bounds on it and its inverse. In particular, we see that, for $\sigma$ sufficiently small,
$$
E_\sigma\equiv \int_{V_I^\sigma\cap\p\B^5_{\delta_0/2}(a_I)}\lf|\nabla_{tan}^2v_\sigma\rg|^2\,d\HF^4\ \leq\ 2\int_{\B^4_\sigma}\lf|\nabla^2\lf(u\circ A_{I,\delta_0/2}\rg)_{\sigma,\rho_0}\rg|^2\,dy\quad<\quad 2 + 2c_{\sS\HH}\ .
$$
By our initial choice (4) of $\delta_0$ we find that, for such $\sigma$,
$$
\int_{ V_I^\sigma\cap\Pi_I^{-1}(0,\delta_0/2]}\lf|\nabla^2v_\sigma\rg|^2\,dx\quad \leq\quad (\delta_0/2)E_\sigma\quad \leq\quad \lf(1+\mathrm{card}(\mathrm{Sing}\,u)\rg)^{-1}\varepsilon/5\ .
$$
In case $b_I\in\mathrm{Sing}\,u$, we make a similar estimate near $b_I$. In any case, we now have 
$$
\limsup_{\sigma\to 0}\sum_{I\in\mathcal{I}}\int_{V_I^\sigma\setminus W_I^\sigma}\lf|\nabla^2 v_\sigma\rg|^2\,dx\ \leq\  \varepsilon/5\ ,
$$
which together with (\ref{transverse}) and (\ref{tintegral}), gives the desired Hessian integral estimate (\ref{wsigest}) for $v_\sigma$. 

Now using  (\ref{wsigest}), we are ready to fix a positive $\sigma_0<1$ so that 
\be\label{wsig0}
\sum_{I\in\mathcal{I}}\int_{V^{\sigma_0}_I}\lf|\nabla^2v_{\sigma_0}\rg|^2\,dx\quad   \leq\quad c_{\sS\HH}\HF^1(\Gamma)\ +\ \varepsilon/2\ .
\ee

The final step will be to modify $v_{\sigma_0}$ to get $u_\varepsilon$. The map $v_{\sigma_0}$ is smooth on $\overline{\B^5}\setminus\mathrm{Sing}\,u$ and is degree-$0$ homogeneous about each point $a\in\mathrm{Sing}\,u$, in the ball $\B^5_{\delta_0/2}(a)$.  

For each such $a$, consider the normalized map given by rescaling $v_{\sigma_0}\,|\,\p\B^5_{\delta_0/2}(a)$, namely,
$$
g_a\ :\ \sS^4\ \to\ \sS^3\ ,\quad g_a(x)\ =\ v_{\sigma_0}[a+(\delta_0/2) x]\ 
$$
We claim that, in $\Pi_4(\sS^3)\simeq \Z_2$,   the homotopy class $\lseg\,g_{a}\rseg$ is zero. 
To see this,
suppose that $a=a_I$ and first note that the restriction of the original map $u\,|\,\p\B^5_{\delta_0/2}(a)$ gives the nonzero class $\lseg\sS\HH\rseg \in \Pi_4(\sS^3)$ by the definition of $\mathcal R$. Second, we slightly reparameterized $u\,|\,\p\B^5_{\delta_0/2}(a)$ near the point $a+(\delta_0/2)\mathbf{e}_I$ to have constant value $\xi_{a_I}=(\sS\HH)(\mathbf{e}_I)$ in a small spherical cap of radius $\sigma_0\delta_0/2$. The resulting reparameterized map $\tilde u_a$ still induces the nonzero homotopy class in $\Pi_4(\sS^3)$. Third, in forming the map $v_\sigma\,|\,\p\B^5_{\delta_0/2}(a)$, we inserted a bubble in the small cap of constancy of $\tilde u_a$. This insertion gives the resulting sum in $\Pi_4(\sS^3)\ $: 
$$
\lseg g_a\rseg\ =\ \lseg \tilde u_a\rseg + \lseg\xi_{a_I}\star\sS\HH\circ\Phi_{\rho_0}\rseg\ =\ \lseg\sS\HH\rseg + \lseg\sS\HH\circ\Phi_{\rho_0}\rseg\ =\ 2\lseg\sS\HH\rseg\ =\ 0\ 
$$
by (\ref{homotopy}).  The same is true in case $a$ is a second endpoint $b_I$.

Now, as in the proof of Lemma \ref{Rdensity}, $g_a$ is homotopic to a constant, and we may we may fix a smooth homotopy 
$h_a:[0,1]\times\sS^4\to\sS^3$ so that
$$
h_a(t,y)\ =\ \cases{ g_a(y) &for $t$ near\ $0$\cr
(1,0,0,0) &for $t$ near\ $1$\cr\ .}
$$
Thus the map 
$$
H_a:\overline{\B^5}\to\sS^3\ ,\quad\ H_a(x)\ =\ h_a\lf(1-|x| ,\,x/|x|\rg)\ \mathrm{for}\ 0<|x|\leq 1\ ,\quad H_a(0)=(1,0,0,0)\ ,
$$
is smooth. Moreover, for $0<\tau\leq\delta_0/2$, 
$$
w_\tau\ :\bigcup_{a\in\mathrm{Sing}\,u}\B^5_\tau(a)\to\sS^3\ ,\quad\ w_\tau(x) = H_a\lf(\frac{x-a}\tau\rg)\ \mathrm{for}\ x\in\B^5_\tau(a)\ ,
$$
satisfies
$$ 
\int_{\B^5_\tau(a)}|\nabla^2w_\tau|^2\,dx\ =\ \tau\int_{\B^5}|\nabla H_a|^2\,dx\ ,
$$
and we can fix a positive $\ \tau_0\leq \delta_0/2\ $ so that
\be\label{tau0}
\sum_{a\in\mathrm{Sing}\,u}\int_{\B^5_{\tau_0}(a)}|\nabla^2w_{\tau_0}|^2\,dx\ <\ \varepsilon/2\ .
\ee

Finally we define the desired map $u_\varepsilon:\B^5\to\sS^3$ by :
$$
u_\varepsilon(x)\ =\ \cases{ v_{\sigma_0}(x) & for $x\in\cup_{I\in\mathcal{I}}V_I^{\sigma_0}\setminus\cup_{a\in\mathrm{Sing}\,u}\B^5_{\tau_0}(a)$\cr
w_{\tau_0}(x) & for $x\in\cup_{a\in\mathrm{Sing}\,u}\B^5_{\tau_0}(a)$\cr
u(x) & otherwise\ .\cr}
$$
We easily verify that $u_\varepsilon$ is smooth and coincides with $u$ outside an $\varepsilon$ neighborhood of $\Gamma$ because $\sigma_0\delta_0/9<\varepsilon$ and
$\tau_0\leq\delta_0/2< \varepsilon$. Moreover, by (\ref{wsig0}) and (\ref{tau0}), 
$$
\begin{array}{lcl}
\ds\int_{\B^5}|\nabla^2u_\varepsilon|^2\,dx\ &\leq& \ds \sum_{I\in\mathcal{I}}\int_{V_I^{\sigma_0}}|\nabla^2v_{\sigma_0}|^2\,dx\ +\sum_{a\in\mathrm{Sing}\,u}\int_{\B^5_{\tau_0}(a)}|\nabla^2w_{\tau_0}|^2\,dx\ +\ \int_{\B^5}|\nabla^2u|^2\,dx\\[5mm]
&\leq&\ds c_{\sS\HH}\HF^1(\Gamma) +\frac \varepsilon 2+\frac \varepsilon 2+  \int_{\B^5}|\nabla^2u|^2\,dx \ .
\end{array}
$$
\cqfd
\begin{Co} If $u$ and $u_\varepsilon$ are as in Theorem \ref{ML}, then $\ u_\varepsilon $ approaches $ u\ $ , $W^{2,2}$ weakly as $\varepsilon\to 0$.
\end{Co}
\vskip.2cm
{\it Proof.}  One has the strong $L^2$ convergence $\lim_{\varepsilon\to 0} \|u_{\varepsilon}-u\|_{L^2} = 0\ $, because the $u_\varepsilon$ are uniformly bounded (by 1) and approach $u$  pointwise on $\B^5\setminus\mathrm{Sing}\,u$. Moreover, for any sequence $1\geq\varepsilon_i\downarrow 0$, we have by Theorem \ref{ML} and Lemma \ref{Sobnorm}, the bound 
$$
\sup_i\| u_{\varepsilon_i}\|^2_{W^{2,2}}\ <\ c_m\lf(4\ + \int_{\B^5}|\nabla^2u|^2\,dx + c_{\sS\HH}\HF^1(\Gamma)\rg)\   <\ \infty\ .
$$
By the weak*(=weak) compactness of the closed ball in $W^{2,2}(\B^5,\R^\ell)$, the sequence $u_{\varepsilon_i}$ contains a subequence $u_{\varepsilon_{i'}}$ that is $W^{2,2}$ weakly convergent to some $w\in W^{2,2}(\B^5,\R^\ell)$. But,  $w$, being by Rellich's theorem,  the strong $L^2$ limit of the $u_{\varepsilon_{i'}}$, must necessarily be the original map  $u$.  Since any subsequence of $u_\varepsilon$ subconverges to the same limit $u$ and since  the weak* (=weak) $W^{2,2}$ topology on bounded sets is metrizable, the original family $u_\varepsilon$ converges $W^{2,2}$ weakly to $u$.
\cqfd

\begin{Rm}\label{csSHH}  In Theorem \ref{ML}, one may replace $c_{\sS\HH}$ by the optimal constant 
$$
\tilde c_{\sS\HH}\ =\ \inf\lf\{\int_{\sS^4}\lf|\nabla^2_{tan} \omega\rg|^2d\HF^4\ : \omega\in\mathcal{C}^\infty(\sS^4,\sS^3)\ \mathrm{and}\ \lseg\omega\rseg=\lseg\sS\HH\rseg \rg\}\ .
$$ 
\end{Rm}

\noindent Here, for any $\omega$ as above, we can first $W^{2,2}$ strongly approximate $u$ by a map  which equals $\omega(x-a)/|x-a|$ in $\B_{\delta_1}(a)$ for all $a\in\mathrm{Sing}\,u$ and some $0<\delta_1<<\delta_0$. Then we repeat the proofs with $\sS\HH$ replaced by $\omega$.

\section{Connecting Singularities with Controlled Length}\label{length}
\medskip

Suppose $u\in{\mathcal R}$ with
Sing$\, u=\{a_1,a_2,\dots,a_m\}$ as above.  Our goal in this
section is to connect the singular points $a_i$, together in pairs or to $\p\B^5$, by some
union of curves whose total length is bounded by an absolute constant
multiple of the {\it Hessian energy}, that is, 
$$
c\int_{\B^5}|\nabla^2u|^2\,dx\ .
$$
This is therefore a bound on the length of a minimal connection for Sing\,$u$, which will allow us, in Theorem \ref{WD} below, to combine Lemma \ref{Rdensity} and Theorem \ref{ML} to obtain the desired sequential weak density of $\mathcal C^\infty(\B^5,\sS^3)$ in $W^{2,2}(\B^5,\sS^3)$.

Using the surjectivity of the suspension of the Hopf map, we readily verify that each {\it regular} value $p\in\sS^3\setminus\{(-1,0,0,0),\,(1,0,0,0)\}$ of $u$
gives a level surface
$$
\Sigma=u^{-1}\{p\}
$$
which necessarily contains all the singular points $a_i$ of $u$.
Note that $\Sigma=u^{-1}\{p\}$ is smoothly embedded away from the
$a_i$ with standard orientation
$\omega_\Sigma\equiv*u^\#\omega_{\sS^3}/|u^\#\omega_{\sS^3}|$, induced from $u$.
Concerning the behavior near $a_i$, the punctured neighborhood
$$
\Sigma\cap\B_{\delta_0}(a_i)\setminus\{a_i\}
$$
is simply a truncated cone whose boundary
$$
\Gamma_i\ =\ \Sigma\cap\p\B_{\delta_0}(a_i)
$$
is a planar circle in the $3$-sphere
$\p\B_{\delta_0}(a_i)\cap\lf(\{\delta p_0\}\times\R^4\rg)$ where $p=(p_0,p_1,p_2,p_3)$.
\medskip

We will eventually choose the desired ``connecting"
curves all to lie on one such level surface $\Sigma$.
\medskip
\subsection{Estimates for Choosing the Level Surface $\Sigma=u^{-1}\{p\}$}\label{surfacesubsection}

We first recall the 3 Jacobian $J_3u = \|\wedge_3Du\|$ and apply
the coarea formula \cite{Fe}, \S3.2.12 with
$$
g\ =\ \frac{|\nabla u|^4+|\nabla^2u|^2}{J_3u}\ ,
$$
to obtain the relation
\be\label{integraldp}
\int_{\sS^3}\int_{u^{-1}\{p\}}\frac{|\nabla u|^4+|\nabla^2u|^2}{J_3u}\,d\HF^2\,d\HF^3p
\ =\ \int_{\B^5}\lf(|\nabla u|^4+|\nabla^2u|^2\rg)dx\ .
\ee
Moreover, since $\|u\|_{L^\infty}=1$, we also have (see \cite{MR}) the
integral inequality
\be\label{w14<w22}
\int_{\B^5}|\nabla u|^4\leq\ c\int_{\B^5}|\nabla^2u|^2\,dx\ .   
\ee
In case $u$ is constant on $\p\B^5$, we may verify this by
computing
$$
\begin{array}{lcl}
\ds\int_{\B^5}|\nabla u|^4 &= &\ds\int_{\B^5}\lf(\nabla u\cdot\nabla u\rg)|\nabla u|^2\,dx
\\[5mm]
& = & \ds\int_{\B^5}\lf[\mbox{div}\lf(u\nabla u|\nabla u|^2\rg)\ -\ u\cdot
\Delta u|\nabla u|^2\ -\
u\nabla u\cdot\nabla\lf(|\nabla u|^2\rg)\rg]\,dx
\\[5mm]
& \leq & \ds 0\ +\ 5\int_{\B^5}|\nabla^2u||\nabla u|^2\,dx\ +\ 2\int_{\B^5}|\nabla^2u||\nabla u|^2\,dx
\\[5mm]
& \leq & \ds\frac 12\int_{\B^5}|\nabla u|^4\,dx\ +\ \frac{49}2\int_{\B^5}|\nabla^2u|^2\,dx .
\end{array}
$$
In the general case, we write $u=\sum_{i=1}^\infty\lambda_iu$ where
$\{\lambda_i\}$ is a partition of unity adapted to a family of Whitney
cubes for $\B^5$.  See \cite{MR}. (The above inequality is true even
with the constraint $\|u\|_{BMO}\leq1$ in place of
$\|u\|_{L^\infty}\leq1$ \cite{MR}.)

By (\ref{integraldp}) and (\ref{w14<w22})  we may now choose a regular value
$p\in\sS^3$ of $u$ so that 
\be\label{regularvaluep}
\int_{u^{-1}\{p\}}\frac{|\nabla u|^4+|\nabla^2u|^2}{J_3u}\,d\HF^2\ \leq\
c\int_{\B^5}|\nabla^2u|^2\,dx\ .                       
\ee
By increasing $c$ we will also  insist that $|p_0|$ is small, say, $\ |p_0|\,<\, 1/100\ $.  This smallness will be useful in guaranteeing that each tangent plane $\mbox{Tan}\,(\Sigma,x)$, for $x\in\Sigma\cap\cup_{i=1}^m\B_{\delta_0}(a_i)\setminus\{ a_i\}$, is close to $\{0\}\times\R^4$.

\medskip
\subsection{A Pull-back Normal Framing for $\Sigma=u^{-1}\{p\}$}.\label{pullback}

Suppose again that
$p=(p_0,p_1,p_2,p_3)\in\sS^3\setminus\{(-1,0,0,0),\,(1,0,0,0)\}$ is a
regular value of $u$. Then
$$
\eta_1\ =\ \lf(-\sqrt{1-p_0^2}\, ,\frac{p_0p_1}{\sqrt{1-p_0^2}}\, ,
\frac{p_0p_2}{\sqrt{1-p_0^2}}\, ,\frac{p_0p_3}{\sqrt{1-p_0^2}}\,\rg)
$$
is the unit vector tangent at $p$ to the geodesic that runs from
$(1,0,0,0)$ through $p$ to $(-1,0,0,0)$.  We may choose two other
vectors
$$
\eta_2\, ,\ \eta_3\ \in\ \mbox{Tan}\lf(\{p_0\}\times\sqrt{1-p_0^2}\ \sS^2,p\,\rg)\
\subset\
 \mbox{Tan}(\sS^3,p)
 $$
so that $\eta_1,\eta_2,\eta_3$ becomes an orthonormal basis for
$\mbox{Tan}(\sS^3,p)$. Since $p$ is a regular value for $u$, these
three vectors lift to three unique smooth linearly independent
normal vectorfields $\tau_1,\tau_2,\tau_3$ along $\Sigma=u^{-1}\{p\}$.  That
is, at each point $x\in\Sigma$,
$$
\tau_j(x)\,\perp\,\Sigma\ \mbox{at }x\ \ \mbox{and } \ Du(x)\lf[\tau_j(x)\rg]\ =\ \eta_j
$$
for $j=1,2,3$.

Near each singularity $a_i$ the lifted vectorfields
$\tau_1,\tau_2,\tau_3$ are also orthonormal.  In fact, for any point $x\in
\Sigma\cap\B_{\delta_0}(a_i)$, $\frac{x_0-a_{i0}}{|x-a_i|}=p_0$, and
\be\label{tau1}
\tau_1(x)\ =\ \lf(-\sqrt{1-p_0^2}\,
,\frac{p_0}{\sqrt{1-p_0^2}}\frac{x_1-a_{i1}}{|x-a_i|}\, ,
\frac{p_0}{\sqrt{1-p_0^2}}\frac{x_2-a_{i2}}{|x-a_i|}\,
,\frac{p_0}{\sqrt{1-p_0^2}}\frac{x_3-a_{i3}}{|x-a_i|}\,\rg)\ .  
\ee
Also $\tau_1(x),\tau_2(x),\tau_3(x)$ are orthonormal for such $x$
because the Hopf map is horizontally orthogonal and the lifts
$\tau_2(x),\tau_3(x)$ are tangent to the $3$ sphere
$\{p_0\}\times\sqrt{1-p_0^2}\,\sS^3$.

On the remainder of the surface $\Sigma\setminus\cup_{i=1}^m\B_{\delta_0}(a_i)$,
the linearly independent vectorfields $\tau_1,\tau_2,\tau_3$ are not
necessarily orthonormal, and we use their Gram-Schmidt
orthonormalizations
$$
\begin{array}{l}
\ds\tilde\tau_1\ =\ \frac{\tau_1}{|\tau_1|}\ ,
\\[5mm]
\ds\tilde\tau_2\ =\
\frac{\tau_2-(\tilde\tau_1\cdot\tau_2)\tilde\tau_1}{|\tau_2-(\tilde\tau_1\cdot\tau_2)\tilde\tau_1|}
=\
\frac{\tau_2-(\tilde\tau_1\cdot\tau_2)\tilde\tau_1}{|\tilde\tau_1\wedge\tau_2|}\ ,
\\[5mm]
\ds\tilde\tau_3\ =\
\frac{\tau_3-(\tilde\tau_1\cdot\tau_3)\tilde\tau_1-(\tilde\tau_2\cdot\tau_3)\tilde\tau_2}
{|\tau_3-(\tilde\tau_1\cdot\tau_3)\tilde\tau_1-(\tilde\tau_2\cdot\tau_3)\tilde\tau_2|}\
=\
\frac{\tau_3-(\tilde\tau_1\cdot\tau_3)\tilde\tau_1-(\tilde\tau_2\cdot\tau_3)\tilde\tau_2}
{|\tilde\tau_1\wedge\tilde\tau_2\wedge\tau_3|}\  ,
\end{array}
$$
which provide an {\it orthonormal framing} for the unit normal bundle of $\Sigma$.

We need to estimate the total variation of these orthonormalizations. Noting that
$|\nabla\lf(\frac \tau{|\tau|}\rg)|\leq2\frac{|\nabla \tau|}{|\tau|}$ for any differentiable $\tau$,
we see that
$$
\begin{array}{lcl}
\ds
|\nabla\tilde\tau_1| & \leq  &\ds 
2\frac{|\nabla\tau_1|}{|\tau_1|}\ \leq\
2\frac{|\nabla\tau_1||\tau_1||\tau_2||\tau_3|}{|\tau_1||\tau_1\wedge\tau_2\wedge\tau_3|}\ 
=\ 2\frac{|\tau_2||\tau_3||\nabla\tau_1|}{|\tau_1\wedge\tau_2\wedge\tau_3|}\ ,
\\[5mm]
\ds |\nabla\tilde\tau_2| & = &\ds
2\lf[\frac{\tau_2-(\tilde\tau_1\cdot\tau_2)\tilde\tau_1}{|\tilde\tau_1\wedge\tau_2|}\rg]\ 
\leq\ 
2\lf[\frac{2|\nabla\tau_2|+2|\tau_2||\nabla\tilde\tau_1|}{|\tau_1\wedge\tau_2||\tau_1|^{-1}}\rg]\
\\[5mm]
&\leq &\ds  8\lf[\frac{|\tau_1\|\nabla\tau_2|+|\tau_2||\nabla\tau_1|}{|\tau_1\wedge\tau_2|}\cdot
\frac{|\tau_1\wedge\tau_2||\tau_3|}{|\tau_1\wedge\tau_2\wedge\tau_3|}\rg]
\
=\ 8\lf[\frac{|\tau_2||\tau_3||\nabla\tau_1|+|\tau_1||\tau_3||\nabla\tau_2|}{|\tau_1\wedge\tau_2\wedge\tau_3|}\rg]\ ,
\\[5mm]
\ds |\nabla\tilde\tau_3| &\leq &\ds 2\lf[\frac{3|\nabla\tau_3|+2|\tau_3||\nabla\tilde\tau_1|+
2|\tau_3||\nabla\tilde\tau_2|}{|\tilde\tau_1\wedge\tilde\tau_2\wedge\tau_3|}\rg]
\\[5mm]
\ds& \leq &\ds 32\lf[\frac{|\nabla\tau_3|+|\tau_3||\tau_1|^{-1}|\nabla\tau_1|+
|\tau_3|\lf(\frac{|\tau_1\|\nabla\tau_2|+|\tau_2||\nabla\tau_1|}{|\tau_1\wedge\tau_2|}\rg)}
{|\frac{\tau_1}{|\tau_1|}\wedge\lf(\frac{\tau_2}{||\tau_1|^{-1}\tau_1\wedge\tau_2|}\rg)\wedge\tau_3|}\rg]
\\[8mm]
\ds& \leq & \ds 32\lf[\frac{|\tau_1||\tau_2||\nabla\tau_3|+|\tau_2||\tau_3||\nabla\tau_1|+
|\tau_1||\tau_3||\nabla\tau_2|}{|\tau_1\wedge\tau_2\wedge\tau_3|}\rg] \ .
\end{array}
$$
Inasmuch as
$$
|\tau_j|\ \leq\ |\nabla u|\ ,\ \ \  |\nabla\tau_j|\ \leq\ |\nabla^2u|\ ,\ \ \
|\tau_1\wedge\tau_2\wedge\tau_3|\ =\ J_3u\ ,
$$
we deduce the general pointwise estimate
$$
|\nabla\tilde\tau_j|\ \leq\ c\frac{|\nabla u|^2|\nabla^2u|}{|J_3u|}\ \leq\
c\frac{|\nabla u|^4+|\nabla^2u|^2}{|J_3u|}\ ,
$$
which we may integrate using (\ref{regularvaluep}) to obtain the variation
estimate along $\Sigma=u^{-1}\{p\}$,
\be\label{intgradtau}
\int_\Sigma|\nabla\tilde\tau_j|\,d\HF^2\ \leq\ c\int_{\B^5}|\nabla^2u|^2\,dx\ . 
\ee
\medskip
\subsection{Twisting of the Normal Frame\label{Q}
$\tilde\tau_1,\tilde\tau_2,\tilde\tau_3$ About Each Singularity $a_i$}.\label{twisting}

First we recall from\cite{MS}, \S 5-6 that the Grassmannian
$$
\tilde G_2(\R^5)
$$
of {\it oriented} 2 planes through the origin in $\R^5$ is a
compact smooth manifold of dimension $6$.  It may be identified
with the set of simple unit $2$ vectors in $\R^5$,
$$
\{v\wedge w\in \wedge_2\R^5\ :\ v\in\sS^4,\ w\in\sS^4,\ v\cdot w=0\}\ .
$$
We will use the distance $|P-Q|$ on $\tilde G_2(\R^5)$ given by
this embedding into $\wedge_2\R^5\approx\R^{10}$.

For a fixed plane $P\in\tilde G_2(\R^5)$, the set of {\it
nontransverse}\ $2$ planes
$$
{\mathcal Q}_P\ =\ \lf\{Q\in \tilde G_2(\R^5)\ :\ P\cap Q\neq\{0\}\rg\}\
$$
is a (Schubert) subvariety of dimension $1+3=4$ because every
$Q\in {\mathcal Q}_P\setminus\{P\}$ equals $v\wedge w$ for some $w\in \sS^4\cap
P$ and some $v\in\sS^4\cap w^\perp$. These subvarieties are all
orthogonally isomorphic and, in particular, have the same finite
$4$ dimensional Hausdorff measure.  Also
$$
Y_P\ =\ \{Q\in{\mathcal Q}_P\ :\ P^\perp\cap Q\neq\{0\}\}
$$
is a closed subvariety of dimension $3$, and ${\mathcal Q}_P\setminus Y_P$
is a smooth submanifold.

Then, near each singularity $a_i$, the set of $2$ planes
nontransverse to the cone $\Sigma\cap\B_{\delta_0}(a_i)\setminus\{a_i\}$,
$$
W\ =\ \bigcup_{x\in \Sigma\cap\B_{\delta_0}(a_i)\setminus\{a_i\}}{\mathcal Q}_{\mbox{Tan}(\Sigma,x)}\ =\
\bigcup_{x\in \Gamma_i}{\mathcal Q}_{\mbox{Tan}(\Sigma,x)}\ ,
$$
has dimension only $1+4=5<6=\mbox{dim\,}\tilde G_2(\R^5)$.  Note also its location, that $W$ is, by the smallness of $|p_0|$, contained in the tubular neighborhood
$$
V\ \equiv\ \{Q\in\tilde G_2(\R^5)\ :\ \mbox{dist}\lf(Q,\tilde G_2(\{0\}\times\R^4\rg) < 1/50\},
$$
of the $4$ dimensional subgrassmannian $\tilde G_2(\{0\}\times\R^4)$.

We now describe explicitly how the framing
$\tilde\tau_1(x),\tilde\tau_2(x),\tilde\tau_3(x)$ {\it twists once} as
$x$ goes around each circle $\Gamma_i$. The problem is that the
vectors $\tilde\tau_j(x)$ lie in the normal space $\mbox{Nor}(\Sigma,x)$
which also varies with $x$. To measure the rotation of the frame
$\tilde\tau_1(x),\tilde\tau_2(x),\tilde\tau_3(x)$, as $x$ traverses the
circle $\Gamma_i$, it is necessary to use some {\it reference frame}
for $\mbox{Nor}(\Sigma,x)$.

We can induce such a frame from some fixed unit vectors in $\R^5$
as follows:  Consider a fixed $Q\in \tilde G_2(\R^5)\setminus W$, and
suppose $Q=v\wedge w$ with $v,w$ being an orthonormal basis for
$Q$.  For each $x\in \Gamma_i$, the orthogonal projections of $v,w$
onto $\mbox{Nor}(\Sigma,x)$ are linearly independent; let
$\sigma_1(x),\,\sigma_2(x)$ be their Gram-Schmidt orthonormalizations. We
then get $\sigma_3(x)$ by using the map $u$ to pull-back the
orientation of $\sS^3$ to $\mbox{Nor}(\Sigma,x)$ so that the resulting
orienting $3$ vector is $\sigma_1(x)\wedge\sigma_2(x)\wedge\sigma_3(x)$ for a
unique unit vector $\sigma_3(x)\in\mbox{Nor}(\Sigma,x)$ orthogonal to
$\sigma_1(x),\,\sigma_2(x)$.  We view
$$
\sigma_1(x)\, ,\ \sigma_2(x)\, ,\ \sigma_3(x)
$$
as the {\it reference frame} determined by the fixed vectors
$v,w$. For each $x\in \Gamma_i$, there is then a unique rotation
$\gamma(x)\in\sS \oO (3)$ so that
$$
\gamma(x)\lf[\sigma_j(x)\rg]\ =\ \tilde\tau_j(x)\ \ \ \mbox{for }\ j=1,2,3\ .
$$

In the next paragraph we will check that $\gamma:\Gamma_i\to\sS\oO (3)$
is a single geodesic circle in $\sS\oO (3)$. The twisting of the
frame $\tilde\tau_1,\tilde\tau_2,\tilde\tau_3$ around the circle $\Gamma_i$
is reflected in the fact that such a circle induces the nonzero
element in $\Pi_1(\sS\oO (3))\simeq\Z_2$.

In the special case $v=(1,0,0,0)$, the normalized orthogonal
projection of $v$ onto $\mbox{Nor}(\Sigma,x)$ is, by (\ref{refnormframe}), simply
$$
\sigma_1(x)\ =\ \tilde\tau_1(x)\ .
$$
So in this case, each orthogonal matrix $\gamma(x)$ is a rotation
about the first axis, and one checks that, as $x$ traverses the
circle $\Gamma_i$ once, these rotations complete a single geodesic
circle in $\sS\oO (3)$. For another choice of $v$, the geodesic
circle $\gamma:\Gamma_i\to \sS\oO (3)$ involves a circle of rotations
about a different axis combined with a single orthogonal change of
coordinates.

\medskip
\subsection{A Reference Normal Framing for $\Sigma=u^{-1}\{p\}$}.\label{refnormframe}

The above calculations near the $a_i$ suggest comparing on the
{\it whole} surface $u^{-1}\{p\}$ the pull-back normal framing
$\tilde\tau_1(x),\tilde\tau_2(x),\tilde\tau_3(x)$ with some reference
normal framing $\sigma_1(x),\sigma_2(x),\sigma_3(x)$ induced by two fixed
vectors $v,w$. Unfortunately, there may not exist fixed vectors
$v,w$ so that the corresponding reference framing
$\sigma_1,\,\sigma_2,\,\sigma_3$ is  defined {\it everywhere} on $\Sigma$. In this
section we show that any orthonormal basis $v,w$ of almost every
oriented 2 plane $Q\in \tilde G_2(\R^5)$ gives a reference framing
on $\Sigma$ which is well-defined and smooth except at finitely many
{\it discontinuities}
$$
b_1,b_2,\dots,b_n\ . $$ We will then need to connect the original
singularities $a_i$ to the $b_j$ (or to $\p\B^5$) and, in \S \ref{connectingaibj}, choose other
curves to connect the $b_j$ to each other (or to $\p\B^5$), with all curves having
total length bounded by a multiple of $\int_{\B^5}|\nabla^2u|^2\,dx$.

To find a suitable $Q=v\wedge w$, we will first rule out the
exceptional planes that contain some nonzero vector normal to $\Sigma$
at some point $x\in\Sigma$. The really exceptional $2$ planes that lie
completely in some normal space
$$
X\ =\ \cup_{x\in\Sigma}X_x\quad\mbox{where}\quad X_x\ =\ \{Q\in \tilde G_2(\R^5)\ :\
Q\subset \mbox{Nor}(\Sigma,x)\}\ .
$$
Then $X$ has dimension at most $2+2=4<6=\mbox{dim\,}\tilde G_2(\R^5)$ because
$\mbox{dim\,}\Sigma=2$ and $\mbox{dim\,}\tilde G_2(\R^3)=2$. The remaining set of
exceptional planes
$$
Y\ =\ \cup_{x\in\Sigma}Y_x\quad\mbox{where}\quad Y_x\ =\ \{Q\in \tilde G_2(\R^5)\ :\
\mbox{dim}\lf(Q\cap\mbox{Nor}(\Sigma,x)\rg)=1\}\
$$
has dimension at most $2+2+1=5<6=\mbox{dim\,}\tilde G_2(\R^5)$ because
$$
Y_x =\ \{e\wedge w\ :\ e\in\sS^4\cap\mbox{Nor}(\Sigma,x)\ \mbox{and }\
w\in\sS^4\cap\mbox{Tan}(\Sigma,x)\}\ .
$$
In terms of our previous notation, $Y_{\mbox{Tan}(\Sigma,x)}=X_x\cup Y_x$.

Any unit vector $e\notin\mbox{Nor}(\Sigma,x)$ has a nonzero orthogonal
projection
$$
e_T(x)
$$
onto $\mbox{Tan}(\Sigma,x)$.

Normalizing
$$
\tilde{e}_T (x)\ =\ \frac{e_T(x)}{|e_T(x)|}\ ,
$$
we find  a unique unit vector $e_\Sigma(x)\in\mbox{Tan}(\Sigma,x)$ orthogonal
to $e_T(x)$ so that $\tilde{e}_T (x)\wedge e_\Sigma(x)$ is the standard
orientation of $\mbox{Tan}(\Sigma,x)$.  Then
$$
e\cdot e_\Sigma(x)\ =\ \lf(e-e_T(x)\rg)\cdot e_\Sigma(x)\ +\ e_T(x)\cdot e_\Sigma(x)\ =\ 0\
+\ 0
$$
because $e-e_T(x)\in\mbox{Nor}(\Sigma,x)$.  Thus,
$$
\tilde{e}_T (x),\,e_\Sigma(x),\tilde\tau_1(x),\,\tilde\tau_2(x),\,\tilde\tau_3(x),
$$
is an orthonormal basis for $\R^5$.

Away from the $4$ dimensional unit normal bundle
$$
{\mathcal N}_\Sigma\ =\ \{(x,e)\ :\ x\in\Sigma,\ e\in\sS^4\cap\mbox{Nor}(\Sigma,x)\}\ ,
$$
we now define the basic map
$$
\Phi\ :\ (\Sigma\times\sS^4)\setminus{\mathcal N}_\Sigma\ \to \ G_2(\R^5)\ ,\ \
\Phi(x,e)=e\wedge e_\Sigma(x)\ ,
$$
to parameterize the planes {\it nontransverse to $\Sigma$}  in $\tilde
G_2(\R^5)\setminus Y$. Incidentally, these do include the $2$
dimensional family of tangent planes
$$
Z\ =\ \{Q\in \tilde G_2(\R^5)\ :\ Q=\mbox{Tan}(\Sigma,x)\ \mbox{for some}\
x\in\Sigma\}\ .
$$
In terms of the notation at the beginning of this section, for any
$2$ plane $Q\notin Y$,
$$
Q\in {\mathcal Q}_{\mbox{Tan}(\Sigma,x)}\ \iff\ Q\ =\ \Phi(x,e)\
\mbox{for\ some }\ e\in\sS^4\setminus\mbox{Nor}(\Sigma,x)\ .
$$
Note that $\Phi(x,-e)=\Phi(x,e)$, and, in fact,
$$
\Phi(x,e')=\Phi(x,e)\in\tilde G_2(\R^5)\setminus Y\ \ \iff\ \ e'=\pm \,e\ .
$$

It is also easy to describe the behavior of $\Phi$ at the singular
set ${\mathcal N}_{\Sigma}$.  A $2$ plane $Q$ belongs to $Y$, that is,
$Q=v\wedge w$ for some $v\in\mbox{Nor}(\Sigma,x)\cap\sS^4$ and
$w\in\mbox{Tan}(\Sigma,x)\cap\sS^4$, if and only if
$Q=\lim_{n\to \infty}\Phi(x_n,v_n)$ for some sequence
$(x_n,v_n)\in(\Sigma\times\sS^4)\setminus{\mathcal N}_\Sigma$ approaching $(x,v)$.  The
map $\Phi$ essentially ``blows-up" the $4$ dimensional ${\mathcal N}_\Sigma$
to the $5$ dimensional $Y$, and, in particular, {\it any smooth
curve in $\tilde G_2(\R^5)$ transverse to $Y$ lifts by $\Phi$ to a
pair of antipodal curves in $\Sigma\times\sS^4$ extending continuously
transversally across ${\mathcal N}_\Sigma$}.

We now choose and fix $Q\in G^2(\R^5)$ so that
\smallskip
\noindent{\it neither $Q$ nor $-Q$ belong to the $5$ dimensional
exceptional set $X\cup Y\cup Z$ and both are regular values of
$\Phi$}.  We may also insist that $Q$ is close to the $3$ dimensional Schubert cycle 
$$
H\ =\ \{ (1,0,\dots,0)\wedge (0,v_1,v_2,v_3,v_4)\ :\ v_1^2+\cdots +v_4^2=1\,\}\ ,
$$
say $\mbox{dist}\,(Q,H)<1/100$.  This will guarantee that $Q$ is well separated from 
the open region $V$ that contains $W$. 
\smallskip

\noindent Since
$$
\mbox{dim}(\Sigma\times\sS^4)=6=\mbox{dim\,}\tilde G_2(\R^5)\ ,
$$
$\Phi^{-1}\{Q,-Q\}$ is a finite set,  say
$$
\Phi^{-1}\{Q,-Q\}\ =\
\{(b_1,e_1),(b_1,-e_1),(b_2,e_2),(b_2,-e_2),\dots,(b_n,e_n),(b_n,-e_n)\}\ .
$$
We now see that the reference framing $\sigma_1(x),\sigma_2(x),\sigma_3(x)$ of
$\mbox{Nor}(\Sigma,x)$ corresponding to any fixed orthonormal basis $ v,\, w$
of $Q$ fails to exist precisely at the points $b_1,b_2,\dots,b_n$.
As before, we now have the smooth mapping
$$
\gamma\ :\ \Sigma\setminus\{a_1,\dots,a_m,b_1,\dots,b_n\}\ \to \ \sS\oO (3)\ ,
$$
which is defined by the condition $\gamma(x)\lf[\sigma_j(x)\rg]\ =\
\tilde\tau_j(x)$ for $j=1,2,3$ or, in column-vector notation,
$$
\gamma\ =\ \lf[\sigma_1\sigma_2\sigma_3\rg]^{-1}\lf[\tilde\tau_1\tilde\tau_2\tilde\tau_3\rg]\ .
$$
\medskip
\subsection{Asymptotic Behavior of $\gamma$ Near  the Singularities
$a_i$ and $b_j$}.\label{asymbeh}

As discussed in \S3.1, the map $u$, the surface
$\Sigma=u^{-1}\{p\}$, the frames $\tilde\tau_1,\tilde\tau_2,\tilde\tau_3$
and $\sigma_1,\sigma_2,\sigma_3$, and the rotation field $\gamma$ are all
precisely known near a singularity $a_i$ in the cone neighborhood
$\Sigma\cap\B_{\delta_0}(a_i)\setminus\{a_i\}$. In particular, $\gamma$ is homogeneous
of degree $0$ on $\Sigma\cap\B_{\delta_0}(a_i)\setminus\{a_i\}$; on its boundary
$\gamma |\,\Gamma_i$ is a constant-speed geodesic circle.

At each $b_j$, the frame $\tilde\tau_1,\tilde\tau_2,\tilde\tau_3$ is
smooth, but the frame $\sigma_1,\sigma_2,\sigma_3$, and hence the rotation
$\gamma$, has an essential discontinuity. Nevertheless, we may deduce
some of the asymptotic behavior at $b_j$ because $\pm Q$ were
chosen to be regular values of $\Phi$. In fact, we'll verify:

\noindent{\it The tangent map $\gamma_j$ of $\gamma$ at $b_j$,
$$
\gamma_j:\mbox{Tan}(\Sigma,b_j)\cap\B_1(0)\to \sS\oO (3)\ ,\ \ \
\gamma_j(x)=\lim_{r\to 0}\gamma\lf[\exp_{b_j}^\Sigma(rx)\rg]\ ,
$$
exists and is the homogeneous degree $0$ extension of some
reparameterization of a geodesic circle in $\sS\oO (3)$}.

\noindent In particular, for small positive $\delta$,
$\gamma\,|\,\lf(\Sigma\cap\p\B_\delta(b_j)\rg)$ is an embedded circle inducing the
nonzero element of $\Pi_1(\sS\oO (3))\simeq\Z_2$.

To check this, we use, as above, the more convenient orthonormal
basis $\{e_j,e_{j\Sigma}\}$ for $Q$; that is,
$$
e_{j\Sigma}=e_{j\Sigma}(b_j)\in\mbox{Tan}(\Sigma,b_j)\ \ \mbox{and }\ \ Q=e_j\wedge
e_{j\Sigma} =\Phi(b_j,\pm e_j)\ .
$$
Then, for $x\in\Sigma$, let
$$
e_j^N(x)\ ,\ \ \ e_{j\Sigma}^{\ N}(x)
$$
denote the orthogonal projections of the fixed vectors
$e_j,\,e_{j\Sigma}\ $ onto $\mbox{Nor}(\Sigma,x)$, and
$$
\hat e_j^N(x)
$$
denote the cross-product of $e_{j\Sigma}^{\ N}(x)$ and $e_j^N(x)$ in
$\mbox{Nor}(\Sigma,x)$. These three vectorfields are smooth near $b_j$
with
$$
e_j^N(b_j)\ \neq\ 0\ ,\ \ \ e_{j\Sigma}^{\ N}(b_j)\ =\ 0 \ ,\ \ \ \
\hat e_j^N(b_j)\ =\ 0\ .
$$
Here our insistence that $\pm Q\notin Z$ guarantees that $Q$ is
not tangent to $\Sigma$ at $b_j$.  Let $g_j$ denote the orthogonal
projection of $\R^5$ onto the 2 plane
$$
P_j\ =\ \mbox{Nor}(\Sigma,b_j)\cap\lf[e_j^N(b_j)\rg]^\perp\ .
$$
Then $G_j(x)=g_j\circ e_j^N(x)$ defines a smooth map from a $\Sigma$
neighborhood of $b_j$  to $P_j$, which has, by the regularity of
$\Phi$ at $(b_j,\tilde e_j)$, a simple, nondegenerate zero at $b_j$
(of degree $\pm 1$).  It follows that as $x$ circulates
$\Sigma\cap\p\B_\delta(b_j)$ once, for $\delta$ small, $G_j(x)$ and similarly
$g_j\circ \hat e_j^N(x)$, circulate $0$ once in $P_j$. Returning to
the original basis $ v,\,w$ of $Q$, we now check that, as $x$
circulates $\Sigma\cap\p\B_\delta(b_j)$ once, the frame
$\sigma_1(x),\sigma_2(x),\sigma_3(x)$ approximately, and asymptotically as
$\delta\to 0$, rotates once about the vector $e_j^N(b_j)$. Since the
frame $\tilde\tau_1(x),\tilde\tau_2(x),\tilde\tau_3(x)$ is smooth at
$b_j$, we see that the map $\gamma$ has, at $b_j$, a tangent map
$\gamma_j$ as described above.
\medskip

\subsection{Connecting the Singularities $a_i$ to the $b_j$ or to $\p\B^5$}.\label{connectingaibj}

Here we will find curves reaching all the $a_i$ and $b_j$.
Concerning the $a_i$, we recall from \cite{Br},\S III,10] that $\sS{\oO}(3)$ is isometric to $\R\pP^3\simeq \sS^3/\{x\sim-x\}$. Any
geodesic circle $\Gamma$ in $\sS\oO (3)$ generates $\Pi_1(\sS{\oO}(3))\simeq\Z_2$ and lifts to a great circle $\tilde\Gamma$ in
$\sS^3$. The rotations at maximal distance from $\Gamma$ form another
geodesic circle $\Gamma^\perp$ and the nearest point retraction
$$
\rho_\Gamma:\sS\oO (3)\setminus \Gamma\ \to \ \Gamma^\perp
$$
is induced by the standard nearest point retraction
$$
\rho_{\tilde \Gamma}:\sS^3\setminus\tilde \Gamma\ \to \ \tilde \Gamma^\perp\ .
$$
In particular,
\be\label{nearestptest}
|\nabla\rho_\Gamma(\zeta)|\ \leq\ \frac c{\mbox{dist}(\zeta,\Gamma)}\ \ \mbox{for }\  \zeta\in\sS{\oO}(3)\ .                                                
\ee
Any geodesic circle $\Gamma'$ in $\sS\oO (3)$ that does {\it not}
intersect $\Gamma$ is mapped diffeomorphically by $\rho_\Gamma$ onto the
circle $\Gamma^\perp$. We deduce that if $\Gamma$ is chosen to miss the
asymptotic circles
$$
\gamma(\Gamma_i)\ \ \ \mbox{and }\ \ \ \gamma_j\lf(\mbox{Tan}(\Sigma,b_j)\cap\sS^4\rg)
$$
associated with the singularities $a_i$ and $b_j$, then, on $\Sigma$,
the composition $\rho_\Gamma\circ \gamma$ maps every sufficiently small circle
$$
	\Sigma\cap\p\B_\delta(a_i)\ \ \mbox{and }\ \ \Sigma\cap\p\B_\delta(b_j)
$$
diffeomorphically onto the circle $\Gamma^\perp$.

Under the identification of $\sS\oO (3)$ with $\R\pP^3$, ${\sS\oO}(4)$ acts transitively by isometry on
$$
{\mathcal G}\ =\ \{\mbox{geodesic circles } \Gamma\subset\sS\oO (3)\ \}\ .
$$
Then ${\mathcal G}$ is compact and admits a positive invariant measure
$\mu_{\mathcal G}$.  For $\mu_{\mathcal G}$ almost every circle $\Gamma$,
$$
\Gamma\cap \gamma(\Gamma_i) = \emptyset\ \ \mbox{for }\ i=1,\dots,m\ ,\ \ \ \  \Gamma\cap
\gamma_j\lf(\mbox{Tan}(\Sigma,b_j)\cap\sS^4\rg) = \emptyset\ \ \mbox{for }\ j=1,\dots,n\ ,
$$
and $\Gamma$ is transverse to the map $\gamma$.  In particular,
$\gamma^{-1}(\Gamma)$ is a finite subset
$$
\{c_1,c_2,\dots,c_\ell\}
$$
of $\Sigma$. For such a circle $\Gamma$ and any regular value $z\in \Gamma^\perp$
of
$$
\rho_\Gamma\circ \gamma\ :\ \Sigma\setminus\{a_1,\dots,a_m,b_1,\dots,b_n,c_1,\dots,c_\ell\}\
\to \ \Gamma^\perp\ ,
$$
the fiber
$$
A\ =\ (\rho_\Gamma\circ \gamma)^{-1}\{z\}
$$
is a smooth embedded $1$ dimensional submanifold with
$$
(\mbox{Clos}\,A)\setminus A\ \subset\
\{a_1,\dots,a_m,b_1,\dots,b_n,c_1,\dots,c_\ell\}\cup\p\B^5\ .
$$

We also can deduce the local behavior of $A$ near each of the
points $a_i,b_j,c_k$.  From the above description of the
asymptotic behavior of $\gamma$ near $a_i$ and $b_j$, we see that
$$
\B_{\delta_0}(a_i)\cap\mbox{Clos}\, A
$$
is simply a {\it single line segment with one \underbar
{endpoint} $a_i$} while
$$
\B_{\delta}(b_j)\cap\mbox{Clos}\, A
$$
is, for $\delta$ sufficiently small, a {\it single smooth segment with
one \underbar{endpoint}  $b_j$}. On the other hand,
$$
\B_{\delta}(c_k)\cap\mbox{Clos}\, A
$$
is, for $\delta$ sufficiently small, a {\it single smooth segment with
an \underbar{interior} point $c_k$}.  To see this, observe that,
for the lifted map $\rho_{\tilde \Gamma}:\sS^3\setminus\tilde \Gamma\to \tilde
\Gamma^\perp$ and any point $\tilde z\in\tilde \Gamma^\perp$, the fiber
$\rho_{\tilde S}^{-1}\{z\}$ is an open great hemisphere, centered at
$z$, with boundary $\tilde \Gamma$. It follows for the downstairs map
$\rho_\Gamma$ that $E_z=\mbox{Clos}\lf(\rho_\Gamma^{-1}\{z\}\rg)\}$ is a full geodesic
$2$ sphere containing $z$ and the circle $\Gamma$. Since the surface
$\gamma(\Sigma)$ intersects the circle $\Gamma$ transversely at a finite set,
this sphere $E_z$ is also transverse to $\gamma(\Sigma)$ near this set.
Thus, for $\delta$ sufficiently small, $\B_\delta(c_k)\cap\mbox{Clos}\,A$,
being mapped diffeomorphically by $\gamma$ onto the intersection
$E_z\cap \gamma\lf(\Sigma\cap\B_\delta(c_k)\rg)$, is an open smooth segment containing
$c_k$ in its interior.

Combining this boundary behavior with the interior smoothness of
the $1$ manifold $A$, we now conclude that
\medskip

{\it $\B^5\cap\mbox{Clos}\,A$ globally consists of disjoint
smooth segments joining pairs of points from 
$$\{a_1,\dots,a_m,b_1,\dots,b_n\}\cup\p\B^5\ .$$  Moreover, each point $a_i$ or $b_j$ 
is the endpoint of precisely one 
segment.}
\medskip

\subsection{Estimating the Length of the Connecting Set $A$}.

The definition of the $A$ depends on many choices:
\begin{itemize}
\item[(1)] the {\it point} $p\in\sS^3$ near $\{0\}\times\R^4$, which determines the surface
$\Sigma=u^{-1}\{p\}$,

\item[(2)]  the {\it vectors} $\eta_2,\eta_2,\eta_3\in\mbox{Tan}(\sS^3,p)$, which determine
the pull-back normal framing $\tilde\tau_1,\tilde\tau_2,\tilde\tau_3$,

\item[(3)]  the {\it vectors} $ v,\, w\in\sS^4$, which determine the reference normal
framing $\sigma_1,\sigma_2,\sigma_3$ and the rotation field
$\gamma=\lf[\sigma_1\sigma_2\sigma_3\rg]^{-1}\lf[\tilde\tau_1\tilde\tau_2\tilde\tau_3\rg]\
 :\ \Sigma\setminus\{b_1,\dots,b_m\}\to \sS\oO (3)$,

\item[(4)]  the {\it circle} $\Gamma\subset\sS\oO (3)$, which determines the
retraction $\rho_\Gamma:\sS\oO (3)\setminus \Gamma\to \Gamma^\perp$, and

\item[(5)] the {\it point} $z\in \Gamma^\perp$, which finally gives
$A=(\rho_\Gamma\circ \gamma)^{-1}\{z\}$.
\end{itemize}

We need to make suitable choices of these to get the desired
length estimate for $A$. In \S \ref{surfacesubsection} we already used one coarea
formula to choose $p\in\sS^3$ to give the basic estimate (\ref{integraldp})
$$
\int_\Sigma\frac{|\nabla u|^4+|\nabla^2u|^2}{J_3u}\,d\HF^2\ \leq\ c\int_{\B^5}|\nabla^2u|^2\,dx\
,
$$
and the pull-back frame estimate (\ref{intgradtau})
$$
\int_\Sigma|\nabla\tilde\tau_j|\,d\HF^2\ \leq\ c\int_{\B^5}|\nabla^2u|^2\,dx\ ,
$$
{\it independent} of the choice of $\eta_1,\eta_2,\eta_3$, then
followed.  For the choice of $z\in S^\perp$, we want to use another coarea
formula, (\cite{Fe}, \S 3.2.22)
\be\label{intdzest}
\int_{\Gamma^\perp}\HF^1(\rho_\Gamma\circ \gamma)^{-1}\{z\}\,dz\ =\
\int_\Sigma|\nabla(\rho_\Gamma\circ \gamma)|\,d\HF^2\ .                       
\ee
To bound the righthand integral, we first use the chain rule and
(\ref{nearestptest}) for the pointwise estimate
\be\label{pointwiseest}
|\nabla(\rho_\Gamma\circ \gamma)(x)|\ =\ |\nabla(\rho_\Gamma)\lf(\gamma(x)\rg)||\nabla\gamma(x)|\ \leq\
\frac c{\mbox{dist}\lf(\gamma(x),\Gamma\rg)}|\nabla\gamma(x)|\ .              
\ee

Next we observe the finiteness of the integral
$$
C\ =\ \int_{\mathcal G}\frac 1{\mbox{dist}(\zeta,\Gamma)}\,d\mu_{\mathcal G}\Gamma\ <\ \infty\ ,
$$
independent of the point $\zeta\in\sS\oO (3)$. To verify this, we
note that $\mu_{\mathcal G}({\mathcal G})<\infty$ and choose a smooth
coordinate chart for $\sS\oO (3)$ near $\zeta$ that maps $\zeta$ to
$0\in\R^3$ and that transforms circles into affine lines in
$\R^3$. Distances are comparable, and an affine line in
$\R^3\setminus\{0\}$ is described by its nearest point $a$ to the origin
and a direction in the plane $a^\perp$. Since
$$
\mu_{\mathcal G}\{\Gamma\in{\mathcal G}\ :\ \zeta\in \Gamma\}\ =\ 0\ ,
$$
the finiteness of $C$ now follows from the finiteness of the $3$
dimensional integral
$$
\int_{\R^3\cap\B_1}\,|y|^{-1}\,dy\ .
$$

We deduce from Fubini's Theorem, (\ref{intdzest}), and (\ref{pointwiseest}) that
$$
\begin{array}{lcl}
\ds\int_{\mathcal G}\int_{\Gamma^\perp}\HF^1(\rho_\Gamma\circ \gamma)^{-1}\{z\}\,dz\,d\mu_{\mathcal G}\Gamma\ &\leq &\ds c\int_\Sigma|\nabla\gamma(x)|\int_{\mathcal G}\frac 1{\mbox{dist}\lf(\gamma(x),\Gamma\rg)}\,d\mu_{\mathcal G}\Gamma\,d\HF^2x\
\\[5mm]
&\leq & \ds c\,C\int_\Sigma|\nabla\gamma(x)|\,d\HF^2x\ .
\end{array}
$$
Thus there exists a $\Gamma\in{\mathcal G}$ and $z\in\Gamma^\perp$ so that
\be\label{Gammaz}
\HF^1(\rho_\Gamma\circ \gamma)^{-1}\{z\}\ \leq\ c\int_\Sigma|\nabla\gamma(x)|\,d\HF^2x\ .
\ee

To estimate the righthand side, recall the matrix formula
$$
\gamma\ =\
\lf[\sigma_1\sigma_2\sigma_3\rg]^{-1}\lf[\tilde\tau_1\tilde\tau_2\tilde\tau_3\rg]\ .
$$
and use Cramer's rule and the product and quotient rules to deduce
the pointwise bound
\be\label{nablagammaptwise}
|\nabla\gamma(x)|\ \leq\ c\sum_{j=1}^3\lf(\,|\nabla\sigma_j(x)|\ +\
|\nabla\tilde\tau_j(x)|\,\rg)\ .                                
\ee
In light of (\ref{intgradtau}), it remains to bound each term
$\int_\Sigma|\nabla\sigma_j(x)|\,d\HF^2x$ for $j=1,2,3$.

For the first one, note that
\be\label{nablasigma1}
|\nabla\sigma_1|\ =\ |\nabla\lf(\frac{v^N}{|v^N|}\rg)|\ \leq\ 2\frac{|\nabla v^N|}{|v^N|}
\ee
where $v^N(x)$ is the orthogonal projection of $v$ onto the normal
space $\mbox{Nor}(\Sigma,x)$ for each $x\in\Sigma$.  The formula
$$
v^N\ =\ \sum_{j=1}^3(v\cdot\tilde\tau_j)\tilde\tau_j\
$$
and the product rule give the pointwise estimate for the
numerator,
\be\label{nablavN}
|\nabla v^N|\ \leq\ c\sum_{j=1}^3|\nabla\tilde\tau_j|\ ,  
\ee
independent of the choice of $v\in\sS^4$.

To estimate the denominator, we let $v^L$ denote the orthogonal
projection of $v$ to any {\it fixed} $3$ dimensional subspace $L$
of $\R^5$, and observe the finiteness
$$
C_1\ =\ \int_{\sS^4}\frac 1{|v^L|}\,d\HF^4v\ <\ \infty\ ,
$$
independent of $L$.  To verify this, we note that the projection
of $\sS^4$ to $L$ vanishes along a great circle, and, near any
point of this circle, the projection is bilipschitz equivalent to
an orthogonal projection of $\R^4$ to $\R^3$.  So the finiteness
of $C_1$ again follows from the finiteness of the $3$ dimensional
integral $\int_{\R^3\cap\B_1}\,|y|^{-1}\,dy$.

By Fubini's Theorem, (\ref{nablasigma1}), (\ref{nablavN}), and (\ref{intgradtau}),
$$
\begin{array}{lcl}
\ds\int_{\sS^4}\int_\Sigma|\nabla\sigma_1(x)|\,d\HF^2x\,d\HF^4v &\leq & \ds2\int_\Sigma|\nabla
v^N(x)|\int_{\sS^4}\frac 1{|v^N(x)|}\,d\HF^4v\,d\HF^2x
\\[5mm]
&\leq & \ds2C_1\int_{\Sigma}|\nabla v^N(x)|\,d\HF^2x\
\\[5mm]
&\leq & \ds c\sum_{j=1}^3\int_{\Sigma}|\nabla\tilde\tau_j(x)|\,d\HF^2x\
\\[5mm]
&\leq & \ds c\int_{\B^5}|\nabla^2u|^2\,dx\ .
\end{array}
$$
So there exists a $v\in\sS^4$ giving the $\sigma_1$ estimate
\be\label{choiceofv}
\int_\Sigma|\nabla\sigma_1(x)|\,d\HF^2x\ \leq\ c\int_{\B^5}|\nabla^2u|^2\,dx\ . 
\ee

Next we observe that $\sigma_2=\frac{w_2}{|w_2|}$ where $w_2(x)$ is the
orthogonal projection onto the $2$ dimensional subspace
$\mbox{Nor}(\Sigma,x)\cap\sigma_1^\perp$.  We again find
\be\label{nablasigma2}
|\nabla\sigma_2|\ =\ |\nabla\lf(\frac{w_2}{|w_2|}\rg)|\ \leq\ 2\frac{|\nabla w_2|}{|w_2|}\ .
\ee
Now the formula
$$
w_2\ =\ \lf[\sum_{j=1}^3(w\cdot\tilde\tau_j)\tilde\tau_j\rg]\ -\
(w\cdot\sigma_1)\sigma_1\ ,
$$
and the product rule give the pointwise estimate for the
numerator,
\be\label{nablaw2}
|\nabla w_2|\ \leq\ c\lf(|\nabla\sigma_1|\ +\ \sum_{j=1}^3|\nabla\tilde\tau_j|\,\rg)\ ,
\ee
independent of the choice $w\in\sS^4$.

To estimate the denominator, we let $w^M$ denote the orthogonal
projection of $w$ to any {\it fixed} $2$ dimensional subspace $M$
of the hyperplane $v^\perp=\sigma_1^\perp$, and observe the finiteness of
the integral
$$
C_2\ =\ \int_{\sS^4\cap v^\perp}\frac 1{|w^M|}\,d\HF^3w\ <\ \infty\ ,
$$
independent of the choices of $v$ or $M$.  To verify this, we note
that the projection of the $3$ sphere $\sS^4\cap v^\perp$ to $M$
vanishes along a great circle, where it is now bilipschitz
equivalent to an orthogonal projection of $\R^3$ to $\R^2$.  So
the finiteness of $C_2$ this time follows from the finiteness of
the $2$ dimensional integral $\int_{\R^2\cap\B_1}\,|y|^{-1}\,dy$.

By Fubini's Theorem, (\ref{intgradtau}), (\ref{nablavN}), (\ref{choiceofv}), (\ref{nablasigma2}) and (\ref{nablaw2}),
$$
\begin{array}{lcl}
\ds\int_{\sS^4\cap v^\perp}\int_\Sigma|\nabla\sigma_2(x)|\,d\HF^2x\,d\HF^3w\ &\leq &\ds 2\int_\Sigma|\nabla w_2(x)|\int_{\sS^4\cap v^\perp}\frac 1{|w_2(x)|}\,d\HF^3w\,d\HF^2x
\\[5mm]
&\leq &\ds 2C_2\int_\Sigma|\nabla w_2(x)|\,d\HF^2x
\\[5mm]
&\leq &\ds  c\int_\Sigma\lf(|\nabla\sigma_1(x)|\ +\
\sum_{j=1}^3|\nabla\tilde\tau_j(x)|\,\rg)\,d\HF^2x
\\[5mm]
&\leq &\ds  c\int_{\B^5}|\nabla^2u|^2\,dx\ .
\end{array}
$$
So there exists a $w\in\sS^4\cap v^\perp$ giving the $\sigma_2$ estimate
\be\label{existsw}
\int_\Sigma|\nabla\sigma_2(x)|\,d\HF^2x\ \leq\ c\int_{\B^5}|\nabla^2u|^2\,dx\ .
\ee
Finally we may use the product rule and the formula
$$
\begin{array}{lcl}
\sigma_3 & = &\ \lf[(\sigma_1\cdot\tilde\tau_2)(\sigma_2\cdot\tilde\tau_3)-
(\sigma_1\cdot\tilde\tau_3)(\sigma_2\cdot\tilde\tau_2)\rg]\tilde\tau_1
\\[5mm]
&+ & \lf[(\sigma_1\cdot\tilde\tau_3)(\sigma_2\cdot\tilde\tau_1)-
(\sigma_1\cdot\tilde\tau_1)(\sigma_2\cdot\tilde\tau_3)\rg]\tilde\tau_2
\\[5mm]
&+ & \lf[(\sigma_1\cdot\tilde\tau_1)(\sigma_2\cdot\tilde\tau_2)-
(\sigma_1\cdot\tilde\tau_2)(\sigma_2\cdot\tilde\tau_1)\rg]\tilde\tau_3\
\end{array}
$$
along with (\ref{intgradtau}), (\ref{choiceofv}), and (\ref{existsw}) to obtain the $\sigma_3$ estimate
\be\label{sigma3est}
\int_\Sigma|\nabla\sigma_3(x)|\,d\HF^2x\ \leq\ c\int_{\B^5}|\nabla^2u|^2\,dx\ .
\ee

Now we may combine (\ref{Gammaz}), (\ref{nablagammaptwise}), (\ref{intgradtau}), (\ref{choiceofv}), (\ref{existsw}), and (\ref{sigma3est}) to obtain the desired length estimate
\be\label{lengthA}
\HF^1(A)\ =\ \HF^1(\rho_\Gamma\circ \gamma)^{-1}\{z\}\ \leq\
c\int_{\B^5}|\nabla^2u|^2\,dx\ .                            
\ee
\medskip

\subsection{Connecting the Singularities $b_j$ to $b_{j'}$}.\label{3.8}

Although we now have a good description and length estimate for
$A$, we are not done.  The problem is that the set $\mbox{Clos}\,A$
does not necessarily connect each of the original singularities
$a_i$ to another $a_{i'}$ or to $\p\B^5$. The path in $\mbox{Clos}\, A$ starting at
$a_i$ may end at some $b_j$. To complete the connections between
pairs of $a_i$, it will be sufficient to find a {\it different}
union $B$  of curves which connect each frame singularity $b_j$ to $\p\B^5$ or to 
another unique frame singularity $b_{j'}$.  Then adding to
$\mbox{Clos}\, A$ some components of $B$ will give the desired curves
connecting every $a_i$ to a distinct $a_{i'}$ or to $\p\B^5$.  In this section we
will use the map $\Phi$ from \S \ref{refnormframe} to construct this additional
connecting set $B$, and we will, in \S \ref{estimateB},
obtain the required estimate on the length of $B$.

First we recall the description in \cite{MS} of $\tilde G_2(\R^5)$ as a
2 sheeted cover of the Grassmannian of {\it unoriented} 2 planes
in $R^5$. With $Q\in\tilde G_2(\R^5)$ chosen as before in \S \ref{Q},
consider the 5 dimensional Schubert cycle
$$
{\mathcal S}_Q\ =\ \{P\in \tilde G_2(\R^5)\ :\ \mbox{dim}\lf(P\cap Q^\perp\rg)\geq 1\}
$$
and the $4$ dimensional subcycle
$$
	{\mathcal T}_Q\ =\ \{P\in \tilde G_2(\R^5)\ :\ \mbox{dim}\lf(P\cap
Q^\perp\rg)\geq 2\}\ =\ \{P\in \tilde G_2(\R^5)\ :\ P\subset Q^\perp\}\ .
$$ 
 As in \cite{MS}, we see that ${\mathcal S}_Q\setminus{\mathcal T}_Q$ is a smooth embedded open 5 dimensional submanifold of  $\tilde G_2(\R^5)$ and that  
$\tilde G_2(\R^5)\setminus {\mathcal S}_Q$ consists of two open $6$ dimensional antipodal
cells, $D_+$ centered at $Q$ and $D_-$ centered at $-Q$.

Next we will carefully define a (nearest-point) retraction map
$$
\Pi_\Q\,:\,\tilde G_2(\R^5)\setminus\{Q,-Q\}\ \to {\mathcal S}_Q\ .
$$
For $P\in D_+\setminus\{Q\}$, there is a unique vector $v\in P\cap \sS^4$
which is at maximal distance in $P\cap\sS^4$ from $Q\cap\sS^4$ and a
unique vector $w$ in $Q\cap\sS^4$ that is closest to $v$; in
particular, $0<w\cdot v<1$. Choose $A_P\in \mbox{so}(5)$ so that the
corresponding rotation $\mbox{exp}\, A_P\in SO(5)$ maps $w$ to $v$ and
maps $\tilde w$ to $\tilde v$ where $P=v\wedge\tilde v$ and
$Q=w\wedge\tilde w$. Thus $\mbox{exp}\, A_P$ maps $Q$ to $P$, preserving
orientation. Here $(\mbox{exp}\, tA_P)(w)$ defines a geodesic circle in
$\sS^4$, and
$$
t_P\ \equiv\ \inf\{t>0\ :\ w\cdot(\mbox{exp}\, tA_P)(w)=0\ \}\ >\ 1\ .
$$
Then $(\mbox{exp}\, 2t_pA_P)(w)=-w$ and $\mbox{exp}\, 4t_pA_P=\mbox{id}$.  It
follows that, in $\tilde G_2(\R^5)$, as $t$ increases,
$$
(\mbox{exp}\, tA_P)(Q)\in D_+\ \mbox{for }\ 0\leq t_P\ \ \mbox{and }\ \ \
(\mbox{exp}\, tA_P)(Q)\in D_-\ \mbox{for }\ t_P<t\leq 2t_P\ ,
$$
$$
(\mbox{exp}\, 0A_P)(Q)=Q,\ \ (\mbox{exp}\, A_P)(Q)=P,\ \ (\mbox{exp}\, t_pA_P)(Q)\in
{\mathcal S}_Q,\ \ (\mbox{exp}\, 2t_pA_P)(Q)=-Q\ ,
$$
and we let $\Pi_Q(P)=(\mbox{exp}\, t_pA_P)(Q)$.

As $P$ approaches $\p D_+={\mathcal S}_Q$, $t_P\downarrow 1$ and
$|\Pi_Q(P)-P|\to 0$. Thus, let
$$
\Pi_Q(P)\ =\ P\ \ \mbox{for }\ \ P\in {\mathcal S}_Q\ .
$$
Also, let
$$
\Pi_Q(P)\ =\ -\Pi_Q(-P)\ \ \mbox{for }\ \ P\in D_-\setminus\{-Q\}\ .
$$

Next recall that the small tubular neighborhood $V$ of the $4$ dimensional subgrassmannian $\tilde G_2(\{0\}\times\R^4)$ was well-separated from $Q$. It follows that $\Pi_Q(V)$ is in a small neighborhood of the $4$ dimensional cycle 
 $\Pi_Q\lf(\tilde G_2(\{0\}\times\R^4)\rg)$. In particular, the $5$ dimensional measure of  $\Pi_Q(V)$ is small, and one easily 
finds $P\in{\mathcal S}_Q\setminus \Pi_Q(V)$ so that $p\notin \Pi_Q(W)$.

For $P\in{\mathcal S}_Q\setminus{\mathcal T}_Q$, the intersection $P\cap
Q^\perp\cap\sS^4$ consists of $2$ antipodal points in $P\cap\sS^4$ that are
uniquely of maximal distance from $Q\cap\sS^4$, and one sees that
$$
\mbox{Clos\,}\Pi_Q^{-1}\{P\}
$$
contains a single semi-circular geodesic arc joining $Q$ and $-Q$. For
almost all $P\in{\mathcal S}_Q\setminus{\mathcal T}_Q$, this semi-circle meets
transversely both $Y$, and, near $\pm Q$, each small surface
$$
\Phi\lf(\,[\Sigma\cap\B_\delta(b_j)]\times\{e_j\}\rg)\ .
$$

We will choose $P\in{\mathcal S}_Q\setminus{\mathcal T}_Q$ also to be a regular
value of $\Pi_Q\circ \Phi$.  Since, near $P$, ${\mathcal S_Q}$ is is a smooth transverse (in fact, 
orthogonal) to $\Pi_Q^{-1}\{P\}$, we find, using  \ref{refnormframe}, that the set
$$
\lf(\Pi_Q\circ \Phi\rg)^{-1}\{P\}\ =\ \Phi^{-1}(\Pi_Q^{-1}\{P\})
$$
is an embedded 1 dimensional submanifold, containing
$\{(b_1,\pm e_1),\dots,(b_m,\pm e_m)\}$.  In small neighborhoods of
any two points $(b_j,e_j),\ (b_j,-e_j)$ the set
$\mbox{Clos}\,(\Pi_Q\circ \Phi)^{-1}\{P\}$ consists of two smooth segments
(antipodal in the $\sS^4$ factor) which both project, under the
projection
$$
p_\Sigma\ :\ \Sigma\times\sS^4\ \to \ \Sigma\ ,
$$
onto a {\it single}  segment in $\Sigma$ which contains $b_j$.
Continuing these two antipodal segments one direction in $(\Pi_Q\circ \Phi)^{-1}\{P\}$ gives antipodal paths whose final endpoints are
$(b_{j'},e_{j'}),\ (b_{j'},-e_{j'})$ for some $j'$ {\it distinct
from} $j$. Here
$$
e_{j'}\wedge e_{j'\Sigma}\ =\ \Phi\lf(b_{j'},\pm e_{j'}\rg)\ =\ -
\Phi\lf(b_j,\pm e_j\rg)\ =\ -e_j\wedge e_{j\Sigma}\ .
$$
Composing either antipodal path with the projection $p_\Sigma$ gives the same
path connecting $b_j$ and $b_{j'}$.  Similarlly, by continuing in the other direction and projecting gives  Thus the whole set
$$
B\ =\ p_\Sigma\lf[(\Pi_Q\circ \Phi)^{-1}\{P\}\rg]\
$$
provides the desired connection in $\Sigma$.

Also note that these two paths upstairs have similar orientations
induced as fibers of the map $\Pi_Q\circ \Phi$. That is, in the notation
of slicing currents \cite{Fe}, \S4.3,
\be\label{slicing}
p_{\Sigma\#}\lf\langle\,\lseg\Sigma\times\sS^4\rseg\,,\,\Pi_Q\circ \Phi\,,\,Q\,\,\rg\rangle\ =\
2(\HF^2\res B)\wedge\vec B\ , 
\ee
where $\vec B$ is a unit tangent vectorfield along $B$ (in the
direction running from $b_j$ to $b_{j'}$).
\medskip

\subsection{Estimating the Length of the Connecting Set $B$}.\label{estimateB}

The definition of $B$ depends on the choices of:

(1) the {\it point} $p\in\sS^3$ near $\{0\}\times\R^4$ which gives the surface $\Sigma=u^{-1}\{p\}$
and the map
$$
\Phi\ :\ (\Sigma\times\sS^4)\setminus{\mathcal N}_\Sigma\ \to \ \tilde G_2(\R^5)\ ,\ \
\Phi(x,e)=e\wedge e_\Sigma(x)\ ,
$$

(2) the $2$ {\it plane} $Q\in \tilde G_2(\R^5)$ near $H$ which determines the
retraction $\Pi_\Q$ of $\tilde G_2(\R^5)\setminus\{Q,-Q\}$ onto the 5
dimensional Schubert cycle ${\mathcal S}_Q$, and

(3) the $2$ {\it plane} $P\in {\mathcal S}_Q\setminus \Pi_Q(V)$ which gives $B =
p_\Sigma\lf[(\Pi_Q\circ \Phi)^{-1}\{P\}\rg]$.
\medskip
\noindent Having chosen $p\in\sS^3$ as before to obtain estimate (\ref{intgradtau}), we
need to chose $Q$ and $P$ to get the desired length estimate for
$B$.

Concerning $Q$, we first readily verify that the retraction $\Pi_Q$
is locally Lipschitz in $\tilde G_2(\R^5)\setminus\{Q,-Q\}$ and deduce
the estimate
\be\label{nablaPiQ}
|\nabla\Pi_Q(S)|\ \leq\ \frac c{|S-Q||S+Q|}\ \ \mbox{for }\ \ S\in\tilde
G_2(\R^5)\setminus\{Q,-Q\}\ .                           
\ee
Using (\ref{slicing}) and \cite{Fe}, 4.3.1, we may integrate the slices to find that
$$
\begin{array}{lcl}
\ds\int_{{\mathcal S}_Q\setminus\Pi_Q(V)}p_{\Sigma\#}\langle\,\lseg\Sigma\times\sS^4\rseg,\Pi_Q\circ \Phi,P\,\rangle\,d\HF^5P\ &\leq &\ds
p_{\Sigma\#}\int_{{\mathcal S}_Q}\langle\,\lseg\Sigma\times\sS^4\rseg,\Pi_Q\circ \Phi,P\,\rangle\,d\HF^5P
\\[5mm]
&=\ &\ds p_{\Sigma\#}\lf(\lseg\Sigma\times\sS^4\rseg\res(\Pi_Q\circ \Phi)^\#\omega_{{\mathcal S}_Q}\rg)\ ,
\end{array}
$$
where $\omega_{{\mathcal S}_Q}$ is the volume element of ${\mathcal S}_Q$. By
(\ref{slicing}) and Fatou's Lemma,
\be\label{lengthBint}
\begin{array}{lcl}
\ds\int_{{\mathcal S}_Q}2\HF^1\lf(p_\Sigma[(\Pi_Q\circ \Phi)^{-1}\{P\}]\rg)\,d\HF^5P&=&\ds
\int_{{\mathcal S}_Q}\M
[p_{\Sigma\#}\langle\lseg\Sigma\times\sS^4]],\Pi_Q\circ \Phi,P\,\rangle]d\HF^5P
\\[5mm]
&\ds \leq&\ds
\M\lf[p_{\Sigma\#}\lf(\,[[\Sigma\times\sS^4\rseg\res(\Pi_Q\circ \Phi)^\#\omega_{{\mathcal S}_Q}\rg)\rg]
\\[5mm]
& = &\ds\sup_{\alpha\in{\mathcal D}^1(\Sigma),|\alpha|\leq1}\int_\Sigma\int_{\sS^4}(\Pi_Q\circ \Phi)^\#\omega_{{\mathcal S}_Q}\wedge
p_\Sigma^\#\alpha\ .                               
\end{array}
\ee

To estimate this last double integral, we recall from \S \ref{Q}
that, for each fixed $x\in\Sigma\setminus\{a_1,\dots,a_m\}$,
$$
\Phi(x,\cdot)\ :\ \sS^4\setminus\mbox{Nor}(\Sigma,x)\ \to \ {\mathcal Q}_x\ \equiv\ {\cal
Q}_{\mbox{Tan}(\Sigma,x)}\setminus Y_x
$$
is a the smooth, orientation-preserving, $2$-sheeted cover map.
Each map $\Phi(x,\cdot)$ depends only on $\mbox{Tan}(\Sigma,x)$, and any two
such maps are orthogonally conjugate. We will derive the formula
\be\label{betadef}
\lf[(\Pi_Q\circ \Phi)^\#\omega_{{\mathcal S}_Q}\wedge p_\Sigma^\#\alpha\rg](x,\cdot)\ =\
\beta(x,\cdot)p_\Sigma^\#\omega_{\Sigma}(x)\wedge\Phi(x,\cdot)^\#\omega_{{\mathcal Q}_x}
\ee
where $\omega_\Sigma$ and $\omega_{{\mathcal Q}_x}$ denote the volume elements of
$\Sigma$ and ${\mathcal Q}_x$ and $\beta(x,\cdot)$ is a smooth function on
$\sS^4\setminus\mbox{Nor}(\Sigma,x)$ satisfying
\be\label{betaest}
\ds|\beta(x,e)|\ \leq\
\frac{c}{|\Phi(x,e)-Q|^5|\Phi(x,e)+Q|^5}\sum_{j=1}^3|\nabla\tilde\tau_j(x)|\
 \ \mbox{for }\ e\in\sS^4\ .                           
\ee
Before proving (\ref{betadef}), note that the decomposition on the
righthand side is not necessarily smooth in $x$ since the
different ${\mathcal Q}_x$ may overlap for $x$ near a critical point
of $\Phi(\cdot,e)$ for some $e\in\sS^4$. Nevertheless, the formula does
imply the measurability of $\beta(x,e)$ in $x$, and so may be
integrated over $\Sigma$.

To derive (\ref{betadef}), we first note that, with the factorization
$\Sigma\times\sS^4$, there are only two terms in the  $(p,q)$ decomposition
of the $5$ form,
$$
(\Pi_Q\circ \Phi)^\#\omega_{{\mathcal S}_Q}\ =\ \Omega_{2,3}\ +\ \Omega_{1,4}\ .
$$
Thus,
\be\label{formdecomp}
(\Pi_Q\circ \Phi)^\#\omega_{{\mathcal S}_Q}\wedge p_\Sigma^\#\alpha\ =\ 0\ +\
\Omega_{1,4}\wedge p_\Sigma^\#\alpha                            
\ee
because the term $\Omega_{2,3}\wedge p_\Sigma^\#\alpha$, being of type
$(2+1,3)$, must vanish.

For each $S=\Phi(x,\pm e)\in {\mathcal Q}_x\setminus Y_x$, we also have the
factorization
$$
\mbox{Tan}(\tilde G_2(\R^5),S)\ =\ \mbox{Nor}({\mathcal Q}_x,S)\,\times\,
\mbox{Tan}({\mathcal Q}_x,S)\ .
$$
Let $\mu_1,\mu_2,\mu_3,\mu_4,\nu_1,\nu_2$ be an orthonormal basis of
$\wedge^1\mbox{Tan}(\tilde G_2(\R^5),S)$ so that
$$
\mu_1,\mu_2,\mu_3,\mu_4 \in\wedge^1\mbox{Tan}({\mathcal Q}_x,S)\ ,\
\nu_1,\nu_2\in\wedge^1\mbox{Nor}({\mathcal Q}_x,S)\ , \
\mu_1\wedge\mu_2\wedge\mu_3\wedge\mu_4=\omega_{{\mathcal Q}_x}(S)\ ;
$$
thus, $0=\nu_1(v)=\nu_2(v)=\mu_1(w)=\mu_2(w)=\mu_3(w)=\mu_4(w)$ whenever
$v\in\mbox{Tan}({\mathcal Q}_x,S)$ and $w\in\mbox{Nor}({\mathcal Q}_x,S)$. We may
expand the $5$ covector
$$
\begin{array}{ll}
\Pi_Q^\#(\omega_{{\mathcal S}_Q})(S)\ =&
\lambda_1\,\nu_2\wedge\mu_1\wedge\mu_2\wedge\mu_3\wedge\mu_4+
\lambda_2\,\nu_1\wedge\mu_1\wedge\mu_2\wedge\mu_3\wedge\mu_4+
\lambda_3\,\nu_1\wedge\nu_2\wedge\mu_2\wedge\mu_3\wedge\mu_4
\\[3mm]
&+\lambda_4\,\nu_1\wedge\nu_2\wedge\mu_1\wedge\mu_3\wedge\mu_4+
\lambda_5\,\nu_1\wedge\nu_2\wedge\mu_1\wedge\mu_2\wedge\mu_4+
\lambda_6\,\nu_1\wedge\nu_2\wedge\mu_1\wedge\mu_2\wedge\mu_3
\end{array}
$$
where
\be\label{lambdaiest}
|\lambda_i|\leq\frac c{|S-Q|^5|S+Q|^5}\ ,                   
\ee
by (\ref{nablaPiQ}).  Applying $\Phi^\#$ (that is, $\wedge^1D\Phi(x,e)\,$) to
all covectors and taking the $(1,4)$ component, we find that only
the first two terms survive so that
\be\label{Omega14}
\begin{array}{lcl}
\Omega_{1,4}(x,e) &= &
\lf[\lambda_1\Phi^\#\nu_2+\lambda_2\Phi^\#\nu_1\rg]_{(1,0)}
\wedge\Phi^\#\mu_1\wedge\Phi^\#\mu_2\wedge \Phi^\#\mu_3\wedge\Phi^\#\mu_4\\[5mm]
&= & \lf[\lambda_1\Phi^\#\nu_2+\lambda_2
\Phi^\#\nu_1\rg]_{(1,0)}\wedge\Phi(x,\cdot)^\#\omega_{{\mathcal Q}_x}(S)\ . 
\end{array}
\ee
Being of type $(2,0)$, the $2$ covector
\be\label{betaform}
\lf(\lf[\lambda_1\Phi^\#\nu_2+\lambda_2\Phi^\#\nu_1\rg]_{1,0}\wedge p_\Sigma^\#\alpha\rg)(x,e)\ =\
\beta(x,e) p_\Sigma^\#\omega_\Sigma(x)                                  
\ee
for some scalar $\beta(x,e)$,  and (\ref{formdecomp}), (\ref{Omega14}), and (\ref{betaform}) now give the
desired formula (\ref{betadef}). This formula readily implies the
smoothness of $\beta(x,\cdot)$ on $\sS^4\setminus\mbox{Nor}(\Sigma,x)$.

To verify the bound (\ref{betaest}), observe that
\be\label{backnui}
|\lf[\Phi^\#\nu_i\rg]_{1,0}|\ =\
\sup_{v\in\sS^4\cap\mbox{Tan}(\Sigma,x)}\nu_i[\nabla_v\Phi(x,e)]\ ,        
\ee
where $\nabla_v\Phi(x,e)=D\Phi_{(x,e)}(v,0)\in\mbox{Tan}(\tilde G_2(\R^5),S)$.
For any unit vector $v\in\mbox{Tan}(\Sigma,x)$ and any $w\in\R^5$,
$$
v\wedge w\in \mbox{Tan}({\mathcal Q}_x,S)
$$
because we may assume $w\notin\mbox{Tan}(\Sigma,x)$ and then choose a
curve $y(t)$ in $\sS^4\cap v^\perp\setminus\mbox{Nor}(\Sigma,x)$ with
$y'(0)=w-(w\cdot v)v$, hence,
$$
v\wedge w\ = v\wedge y'(0)\ =\ \frac d{dt}_{t=0}\lf(v\wedge y(t)\rg)\ =\
-\frac d{dt}_{t=0}\Phi\lf(x,y(t)\rg)\ .
$$
Thus, for any $2$ vector $\xi\in\mbox{Nor}({\mathcal Q}_x,S)$, $|\xi|\ =\
|\xi\wedge v|$; in particular, $|\xi| = |\xi\wedge\tilde{e}_T (x)|$,
$|\xi| = |\xi\wedge e_\Sigma(x)|$, and hence,
$$
|\xi|\ = |\,\xi\wedge\lf(\tilde{e}_T (x)\wedge e_\Sigma(x)\rg)\,|\ .
$$
Since $\nu_i\in\wedge^1\mbox{Nor}({\mathcal Q}_x,S)$ and $|\nu_i|=1$, we now
find that
\be\label{nuigradvphi}
\nu_i[\nabla_v\Phi(x,e)]\ = \nu_i\lf[\lf(\nabla_v\Phi(x,e)\rg)_{\mbox{Nor}(\Sigma,x)}\rg]\ \leq
|\nabla_v\Phi(x,e)\wedge\lf(\tilde{e}_T (x)\wedge e_\Sigma(x)\rg)|\ .  
\ee

Moreover,
\be\label{gradvphi}
\begin{array}{lcl}
|\,\nabla_v\Phi\wedge(\tilde{e}_T \wedge e_\Sigma)|\ &\leq & |\,\lf(\nabla_v(e\wedge
e_\Sigma)\rg)\wedge(\tilde{e}_T \wedge e_\Sigma)|\ \leq\ |\,(\nabla_v
e_\Sigma)\wedge(\tilde{e}_T \wedge e_\Sigma)|\
\\[5mm]
&= & |\,e_\Sigma\wedge\nabla_v(\tilde{e}_T \wedge e_\Sigma)|\ \leq\
|\nabla_v(\tilde{e}_T \wedge e_\Sigma)|\
\\[5mm]
&= & |\,\nabla_v\lf(*(\,\tilde\tau_1\wedge\tilde\tau_2\wedge\tilde\tau_3)\rg)\,|\
 \leq\ c\sum_{j=1}^3|\nabla\tilde\tau_j|\ ,                  
\end{array}
\ee
where $*$ is the Hodge $*:\wedge_3\R^5\to \wedge_2\R^5\approx\R^5$
[F,1.7.8]) The desired pointwise bound (\ref{betaest}) now follows by
combining (\ref{lambdaiest}), (\ref{betaform}), (\ref{backnui}), (\ref{nuigradvphi}) and (\ref{gradvphi}).

For each $x\in\Sigma$, the pull-back $\Phi(x,\cdot)^\#\omega_{{\mathcal Q}_x}$ is
point-wise a positive multiple of the volume form of $\sS^4$. So we
may first integrate over $\sS^4$ and use (\ref{betaest}) to see that
\be\label{integbeta}
\begin{array}{lcl}
\ds\int_{\sS^4}\beta(x,\cdot)\Phi(x,\cdot)^\#\omega_{{\mathcal Q}_x}\ &\leq &\ds
\int_{\sS^4}|\beta(x,\cdot)|\Phi(x,\cdot)^\#\omega_{{\mathcal Q}_x}
\\[5mm]
&\leq &\ds
c\lf(\sum_{j=1}^3|\nabla\tilde\tau_j(x)|\,\rg)\int_{\sS^4}\frac{\Phi(x,\cdot)^\#\omega_{{\cal
Q}_x}}{|\Phi(x,\cdot)-Q|^5|\Phi(x,\cdot)+Q|^5}
\\[5mm]
&= &\ds c\lf(\sum_{j=1}^3|\nabla\tilde\tau_j(x)|\,\rg)\int_{{\cal
Q}_x}\frac{\omega_{{\mathcal Q}_x}(S)}{|S-Q|^5|S+Q|^5}
\cr &\leq &\ds
c\lf(\sum_{j=1}^3|\nabla\tilde\tau_j(x)|\,\rg)\int_{{\cal
Q}_x}\frac{d\HF^4S}{|S-Q|^5|S+Q|^5}\ .         
\end{array}
\ee

To handle the denominator, we note that the Grassmannian $\tilde
G_2(\R^5)$ is a $6$ dimensional homogeneous space, and we readily
use local coordinates to verify that
\be\label{finint}
C_3\ =\ \int_{\tilde G_2(\R^5)}\frac 1{|S-Q|^5|S+Q|^5}\,d\HF^6Q\ <\ \infty\ ,    
\ee
independent of $S$.

Now we recall (\ref{lengthBint}) and fix a sequence of $1$ forms $\alpha_i\in{\cal
D}^1(\Sigma)$ with $|\alpha_i|\leq1$ so that
$$
\M\lf[p_{\Sigma\#}\lf(\,\lseg\Sigma\times\sS^4\rseg\res(\Pi_Q\circ \Phi)^\#\omega_{{\mathcal S}_Q}\rg)\rg]\
 =\ \lim_{i\to \infty}\int_\Sigma\int_{\sS^4}(\Pi_Q\circ \Phi)^\#\omega_{{\mathcal S}_Q}\wedge
p_\Sigma^\#\alpha_i\ ,
$$
let $\beta_i$ be the corresponding function from the formula (\ref{betadef}),
and use (\ref{lengthBint}), Fatou's Lemma, (\ref{betadef}), (\ref{integbeta}), Fubini's Theorem,
(\ref{finint}), and (\ref{intgradtau}) to obtain our final integral estimate
$$
\begin{array}{ll}
&\ds \int_{\tilde G_2(\R^5)}\int_{{\mathcal S}_Q}
2\HF^1\lf(p_\Sigma[(\Pi_Q\circ \Phi)^{-1}\{P\}]\rg)\,d\HF^5P\,d\HF^6Q\
\\[5mm]
&\ds\leq\ \int_{\tilde
G_2(\R^5)}\lim_{i\to \infty}\int_\Sigma\int_{\sS^4}(\Pi_Q\circ \Phi)^\#\omega_{{\cal
S}_Q\setminus\Pi_Q(V)}\wedge p_\Sigma^\#\alpha_i
\\[5mm]
&\ds\leq\ \liminf_{i\to \infty}\int_{\tilde G_2(\R^5)}\int_\Sigma\int_{\sS^4}
\beta_i(x,\cdot)\omega_{\Sigma}(x)\wedge\Phi(x,\cdot)^\#\omega_{{\mathcal Q}_x}
\\[5mm]
 &\ds\leq\
c\int_\Sigma\lf(\sum_{j=1}^3|\nabla\tilde\tau_j(x)|\,\rg) \int_{{\mathcal Q}_x}
\int_{\tilde G_2(\R^5)}\frac{1}{|S-Q|^5|S+Q|^5}\,d\HF^6Q\ d\HF^4S\
d\HF^2x
\\[7mm]
&\ds\leq\ c\, C_3\sum_{j=1}^3\int_\Sigma|\nabla\tilde\tau_j(x)|\,d\HF^2x\quad\leq\quad
c\int_{\B^5}|\nabla^2u|^2\,dx\ .
\end{array}
$$
The Schubert cycles ${\mathcal S}_Q$ are all orthogonally equivalent
and have the same positive $5$ dimensional Hausdorff measure. So we
can use the final integral inequality to choose first a $2$ plane
$Q\in \tilde G_2(\R^5)$ and then a $2$ plane $P\in {\mathcal S}_Q$ so
that the corresponding connecting set
$$
B\ = p_\Sigma[(\Pi_Q\circ \Phi)^{-1}\{P\}]
$$
satisfies the desired length estimate
$$\HF^1(B)\ \leq\
c\int_{\B^5}|\nabla^2u|^2\,dx\ .
$$
\cqfd
\medskip
\begin{Th}[Length Bound]\label{ME} For any $u\in\mathcal{R}$, Sing$\,u$ has a $\Z_2$ connection $\Gamma$ satisfying
$$
\HF^1(\Gamma)\ \leq\  c\int_{\B^5}|\nabla^2u|^2\,dx\ ,
$$ 
for some absolute constant $c$.  
\end{Th}
{\it Proof.} To form the connection $\Gamma$, one takes the union of the curves from $A$ and $B$. The behavior of the individual curves near the points $a_i$ and $b_j$ has been discussed in subsection \ref{asymbeh}.   The set  $A\cup B$ will  pass through each point $b_j$, likely having a corner at $b_j$. One easily replaces the corner with an embedded smooth curve near  $b_j$. Also the curves contributing to $A$ and $B$ may cross. In $\B^5$, it is easy to perturb the curves to eliminate such crossings. The result is the desired $\Z_2$ connection $\Gamma$. (Alternately one can observe that $A\cup B$ already defines a one dimensional  integer multiplicity chain modulo 2 (see \cite{Fe}[4.2.26] which has boundary $\sum_{i=1}^m\lseg a_i\rseg$ relative to $\p\B^5$. As mentioned before, the minimal (mass-minimizing) connection will automatcally consist of non-overlapping intervals. Those that reach $\p\B^5$ meet it orthogonally.

\section{Sequential Weak Density of $W^{2,2}(\B^5,\sS^3)$}\label{sec3}

We are now ready to prove:

\begin{Th}\label{WD} Any map $v$ in $W^{2,2}(\B^5,\sS^3)$ may be approximated in the  $W^{2,2}$ weak topology by a sequence of smooth maps.
\end{Th}\label{wd}
{\it Proof.}  First we may, by Lemma \ref{Rdensity}, chose, for each positive integer $i$, a map  
\noindent$u_i\in\mathcal{R}$ so that 
$\|u_i-v\|_{W^{2,2}} < \frac{1}{i}$; in particular,
$$\label{energybd}
\|u_i-v\|_{L^2}\ <\ \frac{1}{i}\quad\mathrm{and}\quad I\ =\ \sup_i\int_{\B^5}|\nabla^2u_i|^2\,dx\ <\ \infty\ .
$$
Applying Lemma \ref{ML} to each $u_i$, we note that, as $\varepsilon\to 0$, the smooth approximates $u_{i,\varepsilon}$ approach $u_i$ pointwise on $\B^5\setminus\mathrm{Sing}\,u_i$. Inasmuch as the $u_{i,\varepsilon}$ are pointwise bounded (by $1$), Lebesgue's theorem implies
$$
 \|u_{i,\varepsilon}-u_i\|_{L^2}\ \to\ 0\quad\mathrm{as}\quad \varepsilon\to  0 \ .
$$
Thus we can choose a positive $\varepsilon_i$, so that the smooth map $w_i=u_{i,\varepsilon_i}$
has $\|w_i-u_i\|_{L^2}< \frac{1}{i}$; in particular, $\|w_i-v\|_{L^2}< \frac{2}{i}$, and the  smooth maps $w_i$ converge to $v$ strongly in $L^2$.  

On the other hand, by Lemma \ref{Sobnorm} , Lemma \ref{ML}, Theorem \ref{ME}, and (\ref{energybd}),
$$
\sup_i\| w_i\|^2_{W^{2,2}}\ <\ 2c_m(1\ +\ I)\ <\ \infty\ .
$$
By the weak*(=weak) compactness of the closed ball in $W^{2,2}(\B^5,\R^\ell)$, the sequence $w_i$ contains a subequence $w_{i'}$ that is $W^{2,2}$ weakly convergent to some $w\in W^{2,2}(\B^5,\R^\ell)$. But,  $w$, being by Rellich's theorem,  the strong $L^2$ limit of the $w_{i'}$, must necessarily be the original map  $v$. 
\cqfd

\subsection{Least Connection Length $L(v)$}
By Lemma \ref{Rdensity}, one may now define, for any Sobolev map $v\in W^{2,2}(\B^5,\sS^3)$, the nonnegative number
$$
L(v)\ =\ \lim_{\varepsilon\to 0}\inf\{ \HF^1(\Gamma)\ :\ \Gamma\ \mathrm{is\ a}\ \Z_2\,\mathrm{connection\ for\ Sing}\,u\ \mathrm{for\ some}\ u\in\mathcal{R}\ \mathrm{with}\ \|u-v\|_{W^{2,2}}<\varepsilon\}\ .
$$
Any $u\in\mathcal{R}$ has a minimal $\Z_2$ connection, and  $L(u)$ is its length. In general:

\begin{Th} For any $v\in W^{2,2}(\B^5,\sS^3)$, $L(v)=0\ \iff\ v$ is the $W^{2,2}$ strong limit of smooth maps. 
\end{Th} 

{\it Proof.} The sufficiency is immediate from the definition of $L(v)$. To prove the necessity, we assume $L(v)=0$. Then we may choose, for each $i$, a map $u_i\in\mathcal{R}$ along with a $\Z_2$ connection $\Gamma_i$ of Sing$\,u_i$ so that $\|u_i-v\|_{W^{2,2}}<1/i$ and $\HF^1(\Gamma_i)<1/i$. As in the previous proof, there is an $\varepsilon_i<1/i$ so that the smooth maps $w_i=u_{i,\varepsilon_i}$ converge strongly in $L^2$ and weakly in $W^{2,2}$ to $v$. The lower-semicontinuity
$$
\int_{\B^5}|\nabla^2v|^2\,dx\ \leq\ \liminf_{i\to\infty}\int_{\B^5}|\nabla^2w_i|^2\,dx
$$
follows. On the other hand, we have from Theorem \ref{ML} the inequality 
$$
\int_{\B^5}|\nabla^2w_i|^2\,dx\ -\ \int_{\B^5}|\nabla^2u_i|^2\,dx\ \leq\ \frac 1i+\frac{c_{\sS\HH}}i\ .
$$
as well as, from Lemma \ref{Rdensity}, the $W^{2,2}$ strong convergence 
$$
\lim_{i\to\infty}\int_{\B^5}|\nabla^2u_i|^2\,dx\ =\ \int_{\B^5}|\nabla^2v|^2\,dx\ ,
$$
which together imply the upper semi-continuity
$$
\int_{\B^5}|\nabla^2v|^2\,dx\ \geq\ \limsup_{i\to\infty}\int_{\B^5}|\nabla^2w_i|^2\,dx .
$$
The convergence of the total Hessian energies of the $w_i$ to that of $v$, along with the $W^{2,2}$ weak convergence, now implies the $W^{2,2}$ strong convergence of the smooth maps $w_i$ to $v$.
\cqfd

\end{document}